\documentclass[9pt,twoside,a4paper]{article}
\usepackage{mathrsfs}
\usepackage{mathrsfs}
\usepackage{mathrsfs}
\usepackage{mathrsfs}
\usepackage{mathrsfs}
\usepackage{mathrsfs}
\usepackage{mathrsfs}
\newcommand{\liuhao}{\fontsize{8pt}{\baselineskip}\selectfont}

\usepackage{amsfonts}
\usepackage{amsmath,amssymb}
\usepackage{mathrsfs}
\usepackage{bm}
\usepackage{hyperref}

\newtheorem{thm}{Theorem}[section]

\newtheorem{lem}[thm]{Lemma}
\newtheorem{prop}[thm]{Proposition}
\newtheorem{defn}[thm]{Definition}

\newtheorem{rem}[thm]{Remark}

\numberwithin{equation}{section}\allowdisplaybreaks

\def\leq{\leqslant}
\def\geq{\geqslant}

\begin{document}

\title{ {\bf \Large Global weak solution to real bi-fluid system with
		magnetic field}}
\author
{
	{ {Lin Ma$^{3}$ } \ \   Boling Guo$^2$   }
	\ \ Jie Shao$^{1, 2}$\footnote{Corresponding
		author; Email address:
		shaojiehn@foxmail.com. }\\
	{\liuhao $^{1}$ Department of Mathematics, Nanjing University of Science and Technology, Nanjing,  P. R. China}\\
	{\liuhao $^{2}$Institute of Applied Physics and Computational Mathematics,  Beijing ,  P. R. China}\\
	{\liuhao $^{3}$Department of Mathematics, Capital Normal University, Beijing, P. R. China}\\
	\date{}
}
\date{}
\maketitle

\begin{minipage}{13.5cm}
\footnotesize \bf Abstract. \rm  \quad We prove the existence of global  weak solutions with finite energy to some two-fluid systems with magnetic field and  the results suit for corresponding  two-fluid systems. The proof method is mainly inspired by Novotn\'y et al. \cite{novotny2020weak} and Vasseur et al. \cite{vasseur2019global}. For $academic~magnetic$ $ bi\text{-}fluid~ system$, we focus on the case of pressure law with ideal gases and  the new ingredient is that we  can make the proof more explicit without using the preposed hypotheses and remove some unnecessary conditions of the main theorem in Novotn\'y et al. \cite{novotny2020weak}. Meanwhile, the same proof method can be used for $real~magnetic$ $ bi\text{-}fluid~ system$  with pressure law proposed in Vasseur et al. \cite{vasseur2019global}. 

\vspace{9pt}



\end{minipage}

\section{Introduction} \label{se1}

~~~~In this paper, we study
\begin{align}
\partial_{t}\left(\alpha \rho_{+}\right)+&\operatorname{div}\left(\alpha\rho_{+}\mathbf{u}\right)=0, \label{1} \\
\partial_{t}\left((1-\alpha) \rho_{-}\right)+&\operatorname{div}\left(\left(1-\alpha) \rho_{-}\right) \mathbf{u}\right)=0,\label{2}\\
\partial_{t}\left(\left(\alpha \rho_{+}+(1-\alpha) \rho_{-}\right) \mathbf{u}\right)+&\operatorname{div}\left(\left(\alpha \rho_{+}+(1-\alpha) \rho_{-}\right) \mathbf{u} \otimes \mathbf{u}\right)+\nabla \mathfrak{P}_{+}\left(\rho_{+}\right)\nonumber\\
=\mu \Delta \mathbf{u}&+(\mu+\lambda) \nabla \operatorname{div} \mathbf{u}+(\nabla \times  \mathbf{H}) \times  \mathbf{H},\label{3}\\
 \mathbf{H}_{t}-\nabla \times(\mathbf{u} &\times  \mathbf{H})=-\nabla \times(\nu \nabla \times \mathbf{H}),~~\operatorname{div}\mathbf{H}=0,\label{4}\\
 &\mathfrak{P}_{+}\left(\rho_{+}\right)=\mathfrak{P}_{-}\left(\rho_{-}\right), \label{5}
\end{align}
where $0\leq \alpha\leq 1$, $\rho_+\geq 0,\rho_-\geq 0$, $\mathbf{u}$, $\mathbf{H}$  are unknown functions and $\mathfrak{P}_{\pm}\left(\rho_{\pm}\right)=\left(\rho_{\pm}\right)^{\gamma^\pm}$ denotes the pressure function.  $\alpha$ is the rate function of total volume for the first species, $\rho_+$, $\rho_-$ denote the density for two species respectively, while  $\mathbf{u}$, $\mathbf{H}$ is the velocity and magnetic field of the mixture respectively. The constants $\mu, \lambda$ are the shear and bulk (average) viscosities of the mixture and satisfy $\mu>0, 2\mu+3\lambda\geq 0$. The above equations are set in the time-space domain $Q_T=(0,T)\times \Omega$, where $\Omega$ is a sufficiently smooth bounded domain in  $\mathbb{R}^{3}$. 

The equations \eqref{1}-\eqref{5} are endowed with boundary condition that
\begin{equation}
	\left.\left(\mathbf{u},\mathbf{H}\right)\right|_{\partial \Omega}=0\quad \text { for } t \geq 0, \label{6}
\end{equation}
and the initial conditions in $\Omega$
\begin{equation} \label{7}
\begin{aligned}
\alpha\rho_{+}|_{(0, x)} &=\alpha_{0} \rho_{+, 0}(x):=\rho_{0}(x), \\
(1-\alpha) \rho_{-}|_{(0, x)} &=\left(1-\alpha_{0}\right) \rho_{-, 0}(x):=n_{0}(x), \\
\left(\alpha \rho_{+}+(1-\alpha) \rho_{-}\right) \mathbf{u}|_{(0, x)} &=\left(\alpha_{0} \rho_{+,0}+\left(1-\alpha_{0}\right) \rho_{-, 0}\right) \mathbf{u}_{0}(x):=\mathbf{m}_{0},\\
\mathbf{H}|_{(0, x)}: &=\mathbf{H}_{0}(x), ~\operatorname{div} \mathbf{H}_{0}=0.
\end{aligned}
\end{equation}

We call system \eqref{1}-\eqref{7}  the $real~magnetic~bi\text{-}fluid~ system$ following the appellation in \cite{novotny2020weak}. To dispose the system without discussing the rate function $\alpha$ directly, we can set $\rho=\alpha \rho_{+}$,  $n=(1-\alpha) \rho_{-}$, then \eqref{1}-\eqref{4} and \eqref{6}-\eqref{7} can be transformed into
\begin{align}
\partial_{t} \rho&+\operatorname{div}(\rho \mathbf{u}) =0, \label{bi1}\\
\partial_{t} n&+\operatorname{div}(n \mathbf{u}) =0,  \label{bi2}\\
\partial_{t}((\rho+n) \mathbf{u})+\operatorname{div}((\rho+n) \mathbf{u} \otimes \mathbf{u})&+\nabla P(\rho, n) =\mu \Delta \mathbf{u}\nonumber\\
&+(\mu+\lambda) \nabla \operatorname{div} \mathbf{u}+(\nabla \times \mathbf{H}) \times \mathbf{H}, \label{bi3}\\
\mathbf{H}_{t}-\nabla \times(\mathbf{u} \times \mathbf{H})=-&\nabla \times(\nu \nabla \times \mathbf{H}),~~~\operatorname{div} \mathbf{H}=0,\label{bi4}\\ 
	\left.\left(\mathbf{u},\mathbf{H}\right)\right|_{\partial \Omega}&=0 \quad \text { for } t \geq 0, ~\operatorname{div} \mathbf{H}_{0}=0,\label{bi6}\\	
	\rho(0, x)=\rho_{0}(x),~ n(0, x)=n_{0}(x),~(\rho+n) &\mathbf{u}(0, x)=\mathbf{m}_{0}(x),\mathbf{H}(0, x) =\mathbf{H}_{0}(x). \label{bi7}
\end{align}
and \eqref{5} transformed into
\begin{align}
	P(\rho, n)=\mathfrak{P}_{+}&\left(\rho_{+}(\rho,n)\right)=\mathfrak{P}_{-}\left(\rho_{-}(\rho,n)\right).\label{bi5}
\end{align}
In this paper, we call \eqref{bi1}-\eqref{bi7} an $academic~magnetic~bi\text{-}fluid~ system$. When there is no magnetic field, system \eqref{bi1}-\eqref{bi7} with pressure law \eqref{bi5} and 
\begin{align}
	P(\rho,n)=\rho^\gamma+n^\alpha \label{bi8}
\end{align} 
 is studied in \cite{novotny2020weak,vasseur2019global} respectively. We remark that the model of \cite{novotny2020weak} and the one of \cite{vasseur2019global} are different and do not contain each other, as the pressure function of these two models are different.  The $\rho,n$ of \cite{novotny2020weak}  can build direct relationship through  identity \eqref{bi5},
 however the pressure law is an implicit function of $\rho,n$. While the pressure function in \cite{vasseur2019global} is expressed explicitly by convex function  $\rho^\gamma, n^\alpha$,  yet the $\rho, n$ don't have any information from each other.

When there is no magnetic field, \eqref{bi1}-\eqref{bi7} is called $academic~bi\text{-}fluid~system$ and \eqref{1}-\eqref{7} is called $real~bi\text{-}fluid~system$
in  Novotn\'y et al. \cite{novotny2020weak}. The $real~bi\text{-}fluid~system$ is one of the multi-fluid models that raised in \cite{bresch2018multifluid}.  Multi-fluid models are realistic and important for many practical application scenarios, such as turbulent mixing in nuclear industry, reactive flows etc. However, there are rare works for the  multi-fluid models. If shear and bulk (average) viscosities $\mu, \lambda$ are constants, Vassuer et al. \cite{vasseur2019global} obtained global weak solution for $academic~bi\text{-}fluid~system$ with pressure law \eqref{bi8}, when $\bm{\gamma>\frac{9}{5}},\bm{\alpha\geq 1}$. Novotn\'y et al. \cite{novotny2020weak} obtained global weak solution for $real~bi\text{-}fluid~system$ with $\bm{\gamma^+\geq \frac{9}{5}},\bm{\gamma^->0}$ and some other restriction on relationship of $\gamma^+, \gamma^-$.  The proof method in \cite{vasseur2019global} and \cite{novotny2020weak} followed the line in Feireisl et al. \cite{feireisl2001existence} and they developed some technique respectively that handle the bi-fluid pressure law \eqref{5} and \eqref{bi5}. Bresch et al. \cite{bresch2019finite} obtained existence of global in time weak solutions to a compressible
two-fluid Stokes system and their method mainly followed  the novel
compactness tool developed by Bresch et al.\cite{bresch2018global}, which is much different from \cite{feireisl2001existence}. For more results of weak solution and physical description of bi-fluid with constant shear and bulk (average) viscosities, see \cite{vasseur2019global,novotny2020weak,bresch2019finite} and the references therein.  For results on weak solution of bi-fluid with constant shear and bulk (average) viscosities  depending on density and existence of strong solution for multi-fulid, we refer to \cite{bresch2010global,bresch2012global,bresch2018multifluid,huanyao2018review} the references therein.

The main purpose of this article is to derive the global weak solution with finite energy of $real~magnetic~bi\text{-}fluid~ system$ \eqref{1}-\eqref{7} and $academic~magnetic~bi\text{-}fluid~ system$ \eqref{bi1}-\eqref{bi7} with pressure law \eqref{bi5} or \eqref{bi8}. Wen and Zhu \cite{wen2018global} has obtained global strong solution for $academic~magnetic~bi\text{-}fluid~ system$ \eqref{bi1}-\eqref{bi7} with pressure law $P(\rho,n)=\rho+An^\gamma$. As far as we know, ours are the first results of global weak solution for $magnetic~bi\text{-}fluid~ systems$. Meanwhile, compared to the work of \cite{novotny2020weak}, we can remove the some unnecessary condition in main theorem in \cite{novotny2020weak} and make the proof more explicit. As the magnetic field don't make a difference,  our proof and results  can improve those of  $real ~bi\text{-}fluid~ system$ in \cite{novotny2020weak} and $academic~bi\text{-}fluid$ $system$ in \cite{vasseur2019global}, see more in Remark \ref{re1} and \ref{re2}.

  The main part of the article is about proof of the global finite weak solution for $academic~magnet$ $\text{-}ic$ $bi\text{-}fluid~ system$ \eqref{bi1}-\eqref{bi7} with pressure law \eqref{bi5}, which we can follows directly to give the corresponding results of $real~magnetic~bi$$\text{-}fluid~ system$, and we can use the same method to derive  similar results of $academic~magnetic~bi\text{-}fluid~ system$  with pressure law \eqref{bi8}. For the treatment of magnetic field, we mainly follow the line of Hu and Wang \cite{hu2010global}.  Similar with the mono-fluid isentropic compressible Navier-Stokes equation in Feireisl et al.\cite{feireisl2001existence},
 the main difficulty of disposing other terms is to prove the convergence of the pressure law function. However, the pressure function \eqref{bi5} involves two density function $\rho, n$ and we can't get the exact form of  the pressure function. We firstly research the properties of pressure function and then transform it into mono-function during the second and third step of convergence, thus we can use the method in \cite{feireisl2001existence}. To be specific, we find in Lemma  \ref{le1} that the pressure function $P(\rho,n)$ and the first order derivative of $P$ can be bounded by exponential function of $\rho$ and $n$. Then, we use auxiliary
 function  $s$ \eqref{s} to transform $P(\rho,n)$ into $P(ns,n)$, where the $s$ is fixed during the procedure convergence. After these two steps, we are able to obtain the strong convergence of $n$  in a similar way of \cite{feireisl2001existence}. These ideas are mainly inspired by  Novotn\'y et al. \cite{novotny2020weak} and Vasseur-Wen-Yu \cite{vasseur2019global} and we can make some improvements with more precise analysis.

To begin with the proof, we need to investigate the energy inequality. For any smooth solution of system \eqref{bi1}-\eqref{bi5}, we have following energy inequality holds for $0\leq t\leq T$ :
\begin{equation}
E(t)+\int_{0}^{t} \int_{\Omega}\left(\mu|D \mathbf{u}|^{2}+(\lambda+\mu)(\operatorname{divu})^{2}+\nu|\nabla \times \mathbf{H}|^{2}\right) \mathrm{d} x \mathrm{d} s \leq E(0)\label{ener1},
\end{equation}
 where

\begin{align}
E(t)=\int_{\Omega}\left(\frac{1}{2} (\rho+n) \mathbf{u}^{2}\right.&\left.+H_{P}(\rho, n)(t,\cdot)+\frac{1}{2}|\mathbf{H}|^{2}\right) \mathrm{d} x\label{ener2},
\end{align}
\begin{align}\label{ener3}
H_{P}(\rho, n)&:=\rho \int_{1}^{\rho} \frac{P\left(\zeta, \zeta \frac{n}{\rho}\right)}{\zeta^{2}} \mathrm{d} \zeta,~\text{if} ~\rho>0,~H_{P}(0, 0)=0,
\end{align}
and
\begin{equation}
E(0)=\int_{\Omega}\left( \frac{\left|\mathbf{m}_{0}\right|^{2}}{2(\rho_{0}+n_0)}+H_{P}(\rho_0, n_0)+\frac{1}{2}\left|\mathbf{H}_{0}\right|^{2}\right) \mathrm{d} x. \label{ener0}
\end{equation}
Indeed, if $\rho, n, \mathbf{u}, \mathbf{H} $ smooth enough, we can utilize \eqref{bi1}-\eqref{bi3} to get
\begin{align}
\partial_{t}(\frac{(\rho+n) |\mathbf{u}|^2}{2})+\operatorname{div}((\rho+n) \mathbf{u}\frac{|\mathbf{u}|^2}{2})&+\mathbf{u}\nabla P(\rho, n) =\mu \mathbf{u}\Delta \mathbf{u}\nonumber\\
&+(\mu+\lambda)\mathbf{u} \nabla \operatorname{div} \mathbf{u}+((\nabla \times \mathbf{H}) \times \mathbf{H})\cdot\mathbf{u}, \label{ener4}
\end{align}
and multiply \eqref{bi1}, \eqref{bi2} with  $\dfrac{\partial H_{P}(\rho, n)}{\partial\rho}$,  $\dfrac{\partial H_{P}(\rho, n)}{\partial n}$ with  respectively to derive that
\begin{align}
\partial_t(H_{P}(\rho, n))+\operatorname{div}(H_{P}(\rho, n)\mathbf{u})+P(\rho, n)\operatorname{div}\mathbf{u}=0, \label{ener5}
\end{align}
where we used identity
\begin{align*}
	P(\rho,n)=	\rho\frac{\partial H_P(\rho,n)}{\partial\rho}+n\frac{\partial H_P(\rho,n)}{\partial n}-H_P(\rho,n).
\end{align*}
Integrating \eqref{ener4}, \eqref{ener5} on $Q_T$, combining them with \eqref{bi4} and using the following identities
$$
\begin{aligned}
\int_{\Omega}((\nabla \times \mathbf{H}) \times \mathbf{H}) \cdot \mathbf{u} \mathrm{d} x&=-\int_{\Omega}\left(\mathbf{H}^{\top} \nabla \mathbf{u} \mathbf{H}+\frac{1}{2} \nabla\left(|\mathbf{H}|^{2}\right) \cdot \mathbf{u}\right) \mathrm{d} x,\\
\int_{\Omega}(\nabla \times(\nu \nabla \times \mathbf{H}))& \cdot \mathbf{H} \mathrm{d} x =\nu \int_{\Omega}|\nabla \times \mathbf{H}|^{2} \mathrm{~d} x, \\
\int_{\Omega}(\nabla \times(\mathbf{u} \times \mathbf{H})) \cdot \mathbf{H} \mathrm{d} x &=\int_{\Omega} \left(\mathbf{ H}^{\top} \nabla \mathbf{u} \mathbf{H}+\frac{1}{2} \nabla\left(|\mathbf{H}|^{2}\right) \cdot \mathbf{u}\right) \mathrm{d} x,
\end{aligned}
$$
we get \eqref{ener1} immediately, see also \cite{hu2010global}. Using the transformation $\rho=\alpha \rho_{+}$, $n=(1-\alpha) \rho_{-}$, $P(\rho, n)=\mathfrak{P}_{+}\left(\rho_{+}(\rho,n)\right)$, and replacing \eqref{ener3} by
\begin{equation}
	H_{\mathfrak{P}_{+}}(\alpha,\rho_+, \rho_-)=\alpha \rho_{+} \int_{1}^{\alpha \rho_{+}} \frac{\mathfrak{P}_{+}(\rho_{+}(\zeta,\zeta\frac{(1-\alpha)\rho_-}{\alpha \rho_{+}}))}{\zeta^{2}} \mathrm{d} \zeta, ~H_{\mathfrak{P}_{+}}(0, 0,0)=0,\label{est0}
\end{equation}
then, \eqref{ener1} is also the energy inequality for system \eqref{1}-\eqref{5}. As usual, we set the following restriction on the initial data
\begin{align}
	&\inf _{x \in \Omega} \rho_{0} \geq 0, \quad \inf _{x \in \Omega} n_{0} \geq 0, \quad \rho_{0} \in L^{\gamma^+}(\Omega), n_{0} \in L^{\gamma^-}(\Omega),\nonumber
\\
\frac{\bf{m}_{0}}{\sqrt{\rho_{0}+n_{0}}}& \in L^{2}(\Omega), \text { where } \frac{\bf{m}_{0}}{\sqrt{\rho_{0}+n_{0}}}=0 \text { on }\left\{x \in \Omega \mid \rho_{0}(x)+n_{0}(x)=0\right\}. \label{initial}
\end{align}

The definition of weak solution for \eqref{bi1}-\eqref{bi5} is given as
\begin{defn} \label{def1}
 We call $(\rho,n,\mathbf{u},\mathbf{H})$ a finite energy weak solution of \eqref{bi1}-\eqref{bi5} for any $0<T<\infty$, if it satisfies 
 \begin{itemize}
   \item $\rho\in L^\infty(0,T;L^{\gamma^+}(\Omega)),n\in L^\infty(0,T;L^{\gamma^-}(\Omega)),\mathbf{u}\in L^2(0,T;H^1_0(\Omega))$, $\sqrt{(\rho+n)}\mathbf{u}\in L^\infty(0,T;$  $ L^2(\Omega))$, $\mathbf{H} \in L^{2}\left(0, T ; H_{0}^{1}(\Omega)\right) \cap C\left(0, T; L_{\text {weak }}^{2}(\Omega)\right)$,
   \item the energy inequality \eqref{ener1} holds in $\mathcal{D}^{\prime}\left(Q_T\right)$,
   \item $(\rho,n,\mathbf{u},\mathbf{H})$ solves the system \eqref{bi1}-\eqref{bi5} in $\mathcal{D}^{\prime}\left(Q_T\right)$,
   \item the equation \eqref{bi1} and \eqref{bi2} are satisfied in the sense of renormalized solutions, i.e.
$$
\partial_{t} b(f)+\operatorname{div}(b(f) u)+\left[b^{\prime}(f) f-b(f)\right] \operatorname{div} u=0
$$
holds in $\mathcal{D}^{\prime}\left(Q_{T}\right)$, for any $b \in C^{1}(\mathbb{R})$ such that $b^{\prime}(z) \equiv 0$ for all $z \in \mathbb{R}$ large enough, where $f=\rho, n$.
 \end{itemize}
\end{defn}
Similar with definition \ref{def1}, the definition of weak solution for \eqref{1}-\eqref{7} is
\begin{defn} \label{def2}
	 We call $(\alpha,\rho_+,\rho_-,\mathbf{u},\mathbf{H})$ a finite energy weak solution of \eqref{1}-\eqref{5} for any $0<T<\infty$, if
	\begin{itemize}
		\item
		 $\rho_+$ is the unique non-negative solution of equation
		\begin{equation}
			\rho_{+} \left(\rho_{+}\right)^{\frac{\gamma^+}{\gamma^-}}-\left(\rho_{+}\right)^{\frac{\gamma^+}{\gamma^-}} \rho-n\rho_{+}=0
		\end{equation}
     for given function $\rho$, $n$ as in definition \ref{def1}
		and, $\alpha$, $\rho_-$  solve 
		$$\alpha\rho_+=\rho,~\rho_-=\left(\rho_{+}\right)^{\frac{\gamma^+}{\gamma^-}}$$
		 respectively.
		\item $\alpha\rho_+\in L^\infty(0,T;L^{\gamma^+}(\Omega)),(1-\alpha)\rho_-\in L^\infty(0,T;L^{\gamma^-}(\Omega)),\mathbf{u}\in L^2(0,T;H^1_0(\Omega))$,\\ $\sqrt{(\alpha\rho_++(1-\alpha)\rho_-)}\mathbf{u}\in L^\infty(0,T;$  $ L^2(\Omega))$, $\mathbf{H} \in L^{2}\left(0, T ; H_{0}^{1}(\Omega)\right) \cap C\left(0, T; L_{\text {weak }}^{2}(\Omega)\right)$,
		\item the energy inequality \eqref{ener1} with \eqref{est0} holds in $\mathcal{D}^{\prime}\left(Q_T\right)$,
		\item $(\alpha,\rho_+,\rho_-,\mathbf{u},\mathbf{H})$ solves the system \eqref{1}-\eqref{5} in $\mathcal{D}^{\prime}\left(Q_T\right)$,
		\item the equation \eqref{1} and \eqref{2} are satisfied in the sense of renormalized solutions, i.e.
		$$
		\partial_{t} b(f)+\operatorname{div}(b(f) u)+\left[b^{\prime}(f) f-b(f)\right] \operatorname{div} u=0
		$$
		holds in $\mathcal{D}^{\prime}\left(Q_{T}\right)$, for any $b \in C^{1}(\mathbb{R})$ such that $b^{\prime}(z) \equiv 0$ for all $z \in \mathbb{R}$ large enough, where $f=\alpha\rho_+, (1-\alpha)\rho_{-}$.
	\end{itemize}
\end{defn}

We introduce the set that deals with the relation between $\rho$ and $n$
 \begin{align} \label{o}
	\mathscr{O}_{c_0}:=\left\{(\rho, n) \in R^{2} \mid n\in[0, \infty),  \frac{1}{c_0} n\leq \rho\leq c_0 n,~\text{ for  } c_0\geq 1\right\},
\end{align}
which is fundamental important in the analysis next pages. 
Then, we can present the main result for $academic~magnetic~bi\text{-}fluid~ system$ \eqref{bi1}-\eqref{bi7} with pressure law \eqref{bi5}.
\begin{thm} \label{thm1}
For $\gamma^+\geq1, \frac{9}{5}\leq\gamma^-$, assume the initial data satisfies \eqref{initial} and  
$$\frac{1}{{c}_{0}} n_{0} \leq \rho_0 \leq c_{0} n_{0} ~\text { on } ~\Omega$$
for $c_0\geq 1$, then for system \eqref{bi1}-\eqref{bi5}, there exists a finite energy solution defined as definition \ref{def1}. Moreover, we have $(\rho,n) \in \mathscr{O}_{c_0}$, for a.e. $(t,x)\in Q_T$.
\end{thm}

\begin{thm} \label{thm2}
	For $\gamma^+\geq1, \frac{9}{5}\leq\gamma^-$, assume the initial data satisfies \eqref{initial} and
\begin{equation}
	\frac{1}{c_{0}} (1-\alpha_0)\rho_{-,0} \leq  \alpha_0\rho_{+,0} \leq c_{0} (1-\alpha_0)\rho_{-,0}~\text { on } ~\Omega \label{condition1}
\end{equation}
	for $c_0\geq 1$, then for system \eqref{1}-\eqref{5}, there exists a finite energy solution defined as definition \ref{def1}. Moreover, we have $(1-\alpha)\rho_{-} \leq c_{0} \alpha\rho_{+}$ for a.e. $(t,x)\in Q_T$.
\end{thm}
\noindent{\bf{Proof of Theorem \ref{thm2}:}}
 
 We can immediately obtain the result, if we have Theorem \ref{thm1}, use the properties in Subsection \ref{s2.1} and notice the equivalency of system \eqref{1}-\eqref{7} and \eqref{bi1}-\eqref{bi5}.
  \begin{flushright}
 	$\square$
 \end{flushright}

For system \eqref{bi1}-\eqref{bi7} with pressure law  \eqref{bi8}, we have similar result analogous with Theorem \ref{thm1}, whose approach can be used to prove Theorem \ref{thm3} that is parallel to Theorem 1.2 of \cite{vasseur2019global}. 

\begin{defn} \label{def3}
	We call $(\rho,n,\mathbf{u},\mathbf{H})$ a finite energy weak solution of system \eqref{bi1}-\eqref{bi7} with pressure law  \eqref{bi8} for any $0<T<\infty$, if it satisfies
	\begin{itemize}
		\item $\rho\in L^\infty(0,T;L^{\gamma}(\Omega)),n\in L^\infty(0,T;L^{\alpha}(\Omega)),\mathbf{u}\in L^2(0,T;H^1_0(\Omega))$, $\sqrt{(\rho+n)}\mathbf{u}\in L^\infty(0,T;$  $ L^2(\Omega))$, $\mathbf{H} \in L^{2}\left(0, T ; H_{0}^{1}(\Omega)\right) \cap C\left(0, T; L_{\text {weak }}^{2}(\Omega)\right)$,
		\item the energy inequality 
		\begin{align*}
		\frac{d}{d t} \int_{\Omega}\left[\frac{(\rho+n)|\mathbf{u}|^{2}}{2}+G_{\gamma}(
		\rho)+\frac{1}{\alpha-1} n^{\alpha}+\frac{1}{2}|\mathbf{H}|^{2}\right] d x\\
		+\int_{\Omega}\left[\mu|\nabla \mathbf{u}|^{2}+(\mu+\lambda)|\operatorname{div} \mathbf{u}|^{2}+\nu|\nabla \times \mathbf{H}|^{2} \right] d x \leq 0
	\end{align*}
		holds in $\mathcal{D}^{\prime}\left(Q_T\right)$, where
		$$
		G_{\gamma}(\rho)=\left\{\begin{array}{l}
			\rho \ln \rho-\rho+1, \text { for } ~\gamma=1, \\
			\dfrac{\rho^{\gamma}}{\gamma-1}, \quad ~~~~\quad \text { for }  \gamma>1.
		\end{array}\right.
		$$ 
		\item $(\rho,n,\mathbf{u},\mathbf{H})$ solves the system \eqref{bi1}-\eqref{bi5} in $\mathcal{D}^{\prime}\left(Q_T\right)$,
		\item the equation \eqref{bi1} and \eqref{bi2} are satisfied in the sense of renormalized solutions, i.e.
		$$
		\partial_{t} b(f)+\operatorname{div}(b(f) u)+\left[b^{\prime}(f) f-b(f)\right] \operatorname{div} u=0
		$$
		holds in $\mathcal{D}^{\prime}\left(Q_{T}\right)$, for any $b \in C^{1}(\mathbb{R})$ such that $b^{\prime}(z) \equiv 0$ for all $z \in \mathbb{R}$ large enough, where $f=\rho, n$.
	\end{itemize}
\end{defn}

\begin{thm} \label{thm3}
	For $\gamma\geq 1, \frac{9}{5}\leq\alpha$, assume the initial data satisfies \eqref{initial} and
$$\frac{1}{{c}_{0}} n_{0} \leq \rho_0 \leq c_{0} n_{0} ~\text { on } ~\Omega$$
	for $c_0\geq 1$, then for system \eqref{bi1}-\eqref{bi7} with pressure law  \eqref{bi8}, there exists a finite energy solution defined as definition \ref{def3}. Moreover, we have $$\frac{1}{c_0} n\leq \rho\leq c_0 n,$$
	for a.e. $(t,x)\in Q_T$.
\end{thm}

\begin{rem} \label{re1}
	Compared with Theorem 2 of \cite{novotny2020weak}, the Theorem \ref{thm2} here removes the unnecessary condition $\bar{\Gamma}<G$ in \cite{novotny2020weak} with
	\begin{align*}
		G:=\max\{\gamma^{+}&+\gamma_{Bog}^{+},\gamma^{-}+\gamma_{Bog}^{-}\},~~\gamma^\pm_{Bog}:= \min \{\frac{2}{3} \gamma^{\pm}-1, \frac{\gamma^\pm}{3}\},
		\\
		\bar{\Gamma}:=&\max \left\{\gamma^{+}-\frac{\gamma^{+}}{\gamma^{-}}+1, \gamma^{-}+\frac{\gamma^{-}}{\gamma^{+}}-1\right\}.
	\end{align*}
	 Moreover, our proof of  deriving the strong convergence of density $n$ is more plain compared with those of \cite{novotny2020weak}. These improvements result from the more exquisite calculation of properties of $P(\rho,n)$ in Subsection \ref{s2.1} and the  advantage of  relation for $\rho,n$, i.e. \eqref{o}.
	
		When there is no magnetic field in system \eqref{1}-\eqref{7} and \eqref{bi1}-\eqref{bi7}, Theorem \ref{thm1}, \ref{thm2} and \ref{thm3} still hold, since the magnetic field doesn't make difference in the  proof. Thus, our results improves those in \cite{novotny2020weak, vasseur2019global} correspondingly. 
\end{rem}

\begin{rem} \label{re2}
	As $n$ and $\rho$ are symmetric to each other, 
	we can obtain the same results  in Theorem \ref{thm1}, \ref{thm2} and \ref{thm3},
	if we exchange $\rho,n$ with each other.
\end{rem}

\begin{rem}
	 We learnded recently Wen \cite{wen2021global} developed the method in Vasseur et al. \cite{vasseur2019global} and  has  obtained  $real~bi\text{-}fluid$ $system$ and  $academic~bi\text{-}fluid$ $system$ with pressure law \eqref{bi8} for $\gamma^+\geq \frac{9}{5}, ~\gamma^-\geq \frac{9}{5}$, while his proof didn't depend on the assumed condition
	 	$$\frac{1}{\mathrm{c}_{0}} \rho_{0} \leq n_0 \leq \mathrm{c}_{0} \rho_{0} ~\text { on } ~\Omega.$$
	 	We remark that our proof is a different approach compared with those in Wen \cite{wen2021global}.
\end{rem}

The rest of the paper is organized as follows. We study $academic~magnetic~bi\text{-}fluid~ system$ \eqref{bi1}-\eqref{bi7} with pressure law \eqref{bi5} in Sections \ref{s2.}-\ref{s5.}. In Section \ref{s2.}, we investigate the properties of pressure function  \eqref{bi5} and  give some useful lemmas and propositions. Faedo-Galerkin approach for system  is given in Section \ref{s3.}. Passing to the limit of $\epsilon\rightarrow 0^{+}$ and $\delta\rightarrow 0^{+}$ is shown in Sections \ref{s4.} and \ref{s5.} respectively. Finally, we  study $academic~magnetic~bi\text{-}fluid~ system$ \eqref{bi1}-\eqref{bi7} with pressure law \eqref{bi8} and prove Theorem \ref{thm3}  in Section \ref{s6.}.

\section{Preliminaries} \label{s2.}

\subsection{ Properties of the pressure function $P(\rho,n)$} \label{s2.1}

Recall $\rho=\alpha\rho_{+}$,  $n=(1-\alpha) \rho_{-}$, $
\mathfrak{P}_{+}\left(\rho_{+}\right)=\mathfrak{P}_{-}\left(\rho_{-}\right)
$. Defining $ \mathfrak{q}=\mathfrak{P}_{-}^{-1} \circ \mathfrak{P}_{+}$, where $\mathfrak{P}_{-}^{-1}(\cdot)=(\cdot)^{\frac{1}{\gamma^-}}$, then one has 
\begin{align}
\rho_{-}=\mathfrak{q}(\rho_{+})=\rho_{+}^{\frac{\gamma^+}{\gamma^-}},~~
\rho_{+}=\mathfrak{q}^{-1}(\rho_{-})=\rho_{-}^{\frac{\gamma^-}{\gamma^+}}. \label{q}
\end{align}
It follows that
\begin{align}
n\rho_{+}=(1-\alpha) \rho_{-}\rho_{+}=(1-\alpha)\mathfrak{q}(\rho_{+})\rho_{+}\label{rho+0}
\end{align}
i.e.
\begin{align}
\rho_{+} \mathfrak{q}\left(\rho_{+}\right)-\mathfrak{q}\left(\rho_{+}\right) \rho-n\rho_{+}=0. \label{rho+}
\end{align}
From \eqref{q}, \eqref{rho+0} or \eqref{rho+}, direct calculation gives
\begin{equation} \label{q1}
\left\{\begin{array}{l}
     \rho_+=\rho, ~~~~~~~~~~~~~~~~~~~~~~~~~~\text{if}~~n=0,\rho>0,\\
	\mathfrak{q}(\rho_{+})=\rho_-=n,~~~~~~~~~~~~~~~\text{if}~~n>0, \rho=0,\\
	\rho_+=\rho_-=0,~~~~~~~~~~~~~~~~\quad\text{if}~~\rho=n=0.
\end{array}\right.
\end{equation}
For instance, if $n=0,\rho=\alpha\rho_+>0$, we have $\mathfrak{q}^{-1}(\rho_{-})=\rho_+>0$, thus $n=(1-\alpha)\rho_-=0$ indicates that $\alpha=1$ and $\rho_+=\rho$. For the remaining case $\rho>0, n>0$, $ \rho=\alpha\rho_+ $ gives $\rho_{+}>0$. Set $f(\rho_+):=\rho_{+} \mathfrak{q}\left(\rho_{+}\right)-\mathfrak{q}\left(\rho_{+}\right) \rho-n\rho_{+}$, 
then
$$
\begin{aligned}
	f'(\rho_+)&=\rho_{+} \mathfrak{q}'\left(\rho_{+}\right)+ \mathfrak{q}\left(\rho_{+}\right)-\mathfrak{q}'\left(\rho_{+}\right)\rho -n\\
	&=\mathfrak{q}'(\rho_+)\rho_+\left(1-\alpha\right)+\left(\mathfrak{q}\left(\rho_{+}\right)-n\right)>0,~~\text{for } \rho_{+}(1-\alpha)>0
\end{aligned}
$$
and \eqref{rho+} admits a unique solution $ \rho_{+}(\rho, n)$ according to implicit function theorem. Consequently, for any $\rho \geq 0, n \geq 0,$ \eqref{rho+} has a unique solution
$$
\left\{\begin{array}{c}
0<\rho_{+}=\rho_{+}(\rho, n) \in[\rho, \infty) \text { if } \rho>0 \text { or } n>0 \\
\rho_{+}(0,0)=0
\end{array}\right\}
$$
such that 
$$
\rho_{+}(\rho, 0)=\rho, \quad \rho_{+}(0, n)=\mathfrak{q}^{-1}(n).
$$

Based on the analysis above, we can research the properties of $P(\rho,n)$.
\begin{lem} \label{le1} Assume $0\leq\rho, n$,   then $\rho_+(\rho,n), P(\rho,n)$ is continuous function of $\rho, n$ and  $\partial_{\rho} \rho_+ (\rho, n)$,  $ \partial_{n} \rho_+(\rho, n),$  $ \partial_\rho P(\rho, n), ~\partial_n P(\rho, n)$ exist. Moreover,
	\\$(i)$~suppose $0<\rho, 0<n$, then
\begin{align}
\max \left\{\rho, \underline{q} \left(\rho+\mathfrak{q}^{-1}(n)\right)\right\} \leq \rho_{+}(\rho, n) &\leq \bar{q}\left(\rho+\mathfrak{q}^{-1}(n)\right),  \label{rho+1}\\
\underline{q}\leq\partial_{\rho} \rho_{+}(\rho, n)\leq \bar{q},& \label{rr}\\
\underline{C}\left(\rho^{\gamma^+}+n^{\gamma^-}\right) \leq P(\rho, n) &\leq \overline{C}\left(\rho^{\gamma^+}+n^{\gamma^-}\right), \label{P}\\
	0\leq\partial_{n} \rho_{+}(\rho, n) \leq C& \left(\rho^{1-\frac{\gamma^+}{\gamma^-}}+n^{\frac{\gamma^-}{\gamma+}-1}\right), \label{rho+2}\\
0\leq\partial_\rho P(\rho, n) \leq C&\left(\rho^{\gamma^+-1}+n^{\gamma^--\frac{\gamma^-}{\gamma^+}}\right),	\label{Pr}\\
0\leq\partial_{n} P(\rho, n)\leq C& \left(\rho^{\gamma^+-\frac{\gamma^+}{\gamma^-}}+
n^{\gamma^--1}\right),
\label{Pn}
\end{align}
 where $\bar{q},\underline{q},\underline{C}, \overline{C}$ are positive constants to be defined later.
  \\$(ii)$~Suppose $ \rho=0$  or $ n=0$, then \eqref{rho+1}-\eqref{P} still hold. Likewise, \eqref{Pr}, \eqref{Pn} hold if we  set extra condition  $ 1\leq \gamma^+,\gamma^-$, and  \eqref{rho+2} holds if we set $1\leq  \frac{\gamma^-}{\gamma+}  $. 
\end{lem}

\noindent{\bf{Proof:}}

\noindent{\bf{ (i) $\bm{\rho>0, n>0}$ }}

As $\rho>0, n>0$, we can apply implicit function theorem and  take derivation on \eqref{rho+}  with respect to $\rho$  to derive
\begin{align}
\partial_{\rho} \rho_{+}(\rho, n)=\frac{\rho_{+} \mathfrak{q}\left(\rho_{+}\right)}{\rho \mathfrak{q}\left(\rho_{+}\right)+\rho_{+} \mathfrak{q}^{\prime}\left(\rho_{+}\right)\left(\rho_{+}-\rho\right)}
=\frac{1}{\alpha+\frac{\gamma^+}{\gamma^-}\left(1-\alpha\right)}.
\label{rho+3}
\end{align}
It follows from $0\leq \alpha\leq1$ that
$$
\begin{aligned}
\frac{\gamma^+}{\gamma^-}\leq\alpha+\frac{\gamma^+}{\gamma^-}\left(1-\alpha\right)
=\frac{\gamma^+}{\gamma^-}+\alpha(1-\frac{\gamma^+}{\gamma^-})\leq 1,~~\text{if } \frac{\gamma^+}{\gamma^-}\leq 1,\\
1\leq\alpha+\frac{\gamma^+}{\gamma^-}\left(1-\alpha\right)
=\frac{\gamma^+}{\gamma^-}+\alpha(1-\frac{\gamma^+}{\gamma^-})\leq \frac{\gamma^+}{\gamma^-},~~\text{if } \frac{\gamma^+}{\gamma^-}> 1.
\end{aligned}
$$
Thus,
\begin{align}
\underline{q}=:\min\{1,\frac{\gamma^-}{\gamma^+}\}\leq\partial_{\rho} \rho_{+}(\rho, n)\leq \max\{1,\frac{\gamma^-}{\gamma^+}\}:=\overline{q}. \label{rho+4}
\end{align}
Integrating \eqref{rho+4}, one has
$$\rho_{+}(0, n)+\underline{q}\rho \leq \rho_{+}(\rho, n)\leq \bar{q}\rho+\rho_{+}(0, n),$$
which combined with $\rho_{+}(0, n)=\mathfrak{q}^{-1}(n)$ gives \eqref{rho+1}.

Recalling
\begin{align}
P(\rho, n)=\mathfrak{P}_{+}\left(\rho_{+}(\rho, n)\right)=\left(\rho_{+}(\rho, n)\right)^{\gamma^+},
\end{align}
 one can obtain from \eqref{rho+1} that
\begin{align}
\underline{C}\left(\rho^{\gamma^+}+\left(\mathfrak{q}^{-1}(n)\right)^{\gamma^+}\right)\leq \mathfrak{P}_{+}\left(\rho_{+}(\rho, n)\right) \leq \overline{C}\left(\rho^{\gamma^+}+\left(\mathfrak{q}^{-1}(n)\right)^{\gamma^+}\right), \label{P+}
\end{align}
where $ \underline{C}=(\underline{q})^{\gamma^+}, \overline{C}=(2\overline{q})^{\gamma^+}  $ and which combined with $\left(\mathfrak{q}^{-1}(n)\right)^{\gamma^+}=n^{\gamma^-} $ gives \eqref{P}.

Applying implicit function theorem and taking derivation on \eqref{rho+} with respect to $n$, one has 
\begin{align}
	\begin{aligned}\label{rho+n}
0<\partial_{n} \rho_{+}(\rho, n)&=\frac{\left(\rho_{+}\right)^{2}}{\rho \mathfrak{q}\left(\rho_{+}\right)+\rho_{+} \mathfrak{q}^{\prime}\left(\rho_{+}\right)\left(\rho_{+}-\rho\right)}\\
&=\frac{\left(\rho_{+}\right)^{2}}{(\rho_{+})^2 \mathfrak{q}^{\prime}\left(\rho_{+}\right)}\frac{1}{\frac{\gamma^-}{\gamma^+}\alpha+(1-\alpha)}\\
&\leq \overline{q}_1\frac{\gamma^-}{\gamma^+}(\rho_+)^{-(\frac{\gamma^+}{\gamma^-}-1)}\\
&\leq  \overline{q}_1\frac{\gamma^-}{\gamma^+}\max\{\underline{q}^{-(\frac{\gamma^+}{\gamma^-}-1)},\bar{q}^{-(\frac{\gamma^+}{\gamma^-}-1)}\}\left(\rho+\mathfrak{q}^{-1}(n)\right)^{-(\frac{\gamma^+}{\gamma^-}-1)}\\
&\leq C\left(\rho^{1-\frac{\gamma^+}{\gamma^-}}+n^{\frac{\gamma^-}{\gamma+}-1}\right),
\end{aligned}
\end{align}
where 
\begin{align} \label{q11}
	\underline{q}_1=:\min\{1,\frac{\gamma^+}{\gamma^-}\}\leq\frac{1}{\frac{\gamma^-}{\gamma^+}\alpha+(1-\alpha)}\leq \max\{1,\frac{\gamma^+}{\gamma^-}\}:=\overline{q}_1
\end{align}
and we have used $\mathfrak{q}(\cdot)=(\cdot)^{\frac{\gamma^+}{\gamma^-}}$, $\mathfrak{q}^{-1}(\cdot)=(\cdot)^{\frac{\gamma^-}{\gamma^+}}$ and \eqref{rho+1}.

In the same manner,  we can use \eqref{rho+1}, \eqref{rho+4} and $\mathfrak{q}^{-1}(\cdot)=(\cdot)^{\frac{\gamma^-}{\gamma^+}}$ to derive 
\begin{align}
	0<\partial_\rho P(\rho, n) &=\mathfrak{P}'_{+}\left(\rho_{+}(\rho, n)\right)\partial_{\rho} \rho_{+}(\rho, n)\nonumber\\
	&\leq C\overline{q}\left(\rho+\mathfrak{q}^{-1}(n)\right)^{\gamma^+-1}\nonumber\\
	&\leq C(\underline{r})\overline{q}\left(\rho^{\gamma^+-1}+n^{\gamma^--\frac{\gamma^-}{\gamma^+}}\right), \label{Pr1}\\
0<\partial_{n} P(\rho, n)&=\mathfrak{P}'_{+}\left(\rho_{+}(\rho, n)\right)\partial_{n} \rho_{+}(\rho, n)\nonumber\\
&\leq \mathfrak{P}'_{+}\left(\rho_{+}(\rho, n)\right)  \frac{1}{ \mathfrak{q}^{\prime}\left(\rho_{+}\right)}\nonumber\\
&\leq C \left(\rho+\mathfrak{q}^{-1}(n)\right)^{\gamma^+-\frac{\gamma^+}{\gamma^-}}\nonumber\\
&\leq C \left(\rho^{\gamma^+-\frac{\gamma^+}{\gamma^-}}+
n^{\gamma^--1}\right),
\end{align}
where we have used $\mathfrak{P}'_{+}(\cdot)=\gamma^+(\cdot)^{\gamma^+-1}$, \eqref{rho+1},  $\mathfrak{q}(\cdot)=(\cdot)^{\frac{\gamma^+}{\gamma^-}}$, $\mathfrak{q}^{-1}(\cdot)=(\cdot)^{\frac{\gamma^-}{\gamma^+}}$.

\noindent{\bf{ (ii) $\bm{\rho=0}$ or $\bm{n=0}$}}

We firstly discuss the case $n>0,\rho=0$, where we have $ \mathfrak{q}(\rho_+)=n, \alpha=0$ by following the analysis in \eqref{q1}. Then, 
\begin{align}
\rho_+=\mathfrak{q}^{-1}(n),~~P(0,n)=\rho_+^{\gamma^+}=n^{\gamma^-},
\end{align}
 therefore  \eqref{rho+1}, \eqref{P} holds and direct calculation shows
 \begin{align}
 	\partial_n \rho_+(0,n)=\frac{\gamma^-n^{\frac{\gamma^-}{\gamma^+}-1}}{\gamma^+},~~ \partial_n P(0,n)=\gamma^- n^{\gamma^--1}.
 \end{align}
By \eqref{q} and the definition of $ \rho, n, \rho_{+}, \mathfrak{q}(\cdot) $, one has
\begin{align*}
		\partial_\rho \rho_+(0,n)&=	\lim _{\delta\rho \rightarrow 0^+}\frac{\rho_+(\delta\rho,n)-\rho_+(0,n)}{\delta\rho}
	=	\lim _{\delta\rho \rightarrow 0^+}\frac{\dfrac{\delta\rho}{\alpha}-\mathfrak{q}^{-1}(n)}{\delta\rho}\\
	&=	\frac{1-(1-\alpha)^\frac{\gamma^-}{\gamma^+}}{\alpha},
\end{align*}
where 
\begin{align*}
\mathfrak{q}^{-1}(n)&=(1-\alpha)^\frac{\gamma^-}{\gamma^+} (\rho_-(\delta\rho,n))^\frac{\gamma^-}{\gamma^+}=(1-\alpha)^\frac{\gamma^-}{\gamma^+} \rho_+(\delta\rho,n)\\
&=\frac{(1-\alpha)^\frac{\gamma^-}{\gamma^+}}{\alpha}\delta\rho.
\end{align*}
It is easy to prove that $ g(\alpha):=\dfrac{1-(1-\alpha)^\frac{\gamma^-}{\gamma^+}}{\alpha}, \alpha\in[0,1] $ is decreasing for $ \frac{\gamma^-}{\gamma^+}>1 $ and increasing for  $ 0<\frac{\gamma^-}{\gamma^+}\leq1 $. We also have
\begin{align*}
\lim _{\alpha \rightarrow 0^+}=	\frac{1-(1-\alpha)^\frac{\gamma^-}{\gamma^+}}{\alpha}=\frac{\gamma^-}{\gamma^+},
\end{align*}
thus
\begin{align*}
	\underline{q}\leq\partial_{\rho} \rho_{+}(0, n)\leq \overline{q}.
\end{align*}
As $\delta\rho=\alpha\rho_{+}(\delta\rho,n)$,   one has 
$$\delta\rho\rightarrow0 \Longleftrightarrow\alpha \rightarrow 0,$$
otherwise, if $ \rho_{+}\rightarrow 0 $, there is $n=(1-\alpha) \rho_{-}= (1-\alpha) \rho_{+}^{\frac{\gamma^+}{\gamma^-}}\rightarrow 0$, which contradicts to $ n>0  $. Thus, one has
\begin{align*}
	\partial_\rho P(0,n)&=	\lim _{\delta\rho \rightarrow 0}\frac{P(\delta\rho,n)-P(0,n)}{\delta\rho}\\
	&=	\lim _{\alpha \rightarrow 0}(1-\alpha)^\frac{\gamma^-}{\gamma^+}\frac{\left(\frac{n}{1-\alpha}\right)^{\gamma^-}-n^{\gamma^-}}{\alpha\mathfrak{q}^{-1}(n)}\\
	&=\lim _{\alpha \rightarrow 0} (1-\alpha)^{\frac{\gamma^-}{\gamma^+}-\gamma^-} \frac{1-(1-\alpha)^{\gamma^-}}{\alpha}\frac{n^{\gamma^-}}{\mathfrak{q}^{-1}(n)}\\
	&=\gamma^- n^{\gamma^--\frac{\gamma^-}{\gamma^+}}.
\end{align*}

Similarly for the case  $n=0,\rho>0$, we have  $ \alpha=1,  \rho_+(\rho,0)=\rho, P(\rho,0)=\rho^{\gamma^+}, 	\partial_\rho \rho_+(\rho,0)=1, \partial_\rho P(\rho,0)=\gamma^+\rho^{\gamma^+-1}   $ and
\begin{align*}
	\partial_n \rho_+(\rho,0)&=	\lim _{\delta n \rightarrow 0^+}\frac{\rho_+(\rho,\delta n)-\rho_+(\rho,0)}{\delta n}\\
&=	\lim _{\alpha \rightarrow 1^-}\frac{(\frac 1\alpha-1)\rho}{(1-\alpha)\rho^{\frac{\gamma^+}{\gamma^-}}}\alpha^\frac{\gamma^+}{\gamma^-}\\
&=\rho^{1-\frac{\gamma^+}{\gamma^-}},\\
	\partial_n P(\rho,0)&=	\lim _{\delta n \rightarrow 0}\frac{P(\rho,\delta n)-P(\rho,0)}{\delta n}\\
&=	\lim _{\alpha \rightarrow 1^-}\frac{(\frac {1}{\alpha^{\gamma^+}}-1)\rho^{\gamma^+}}{(1-\alpha)\rho^{\frac{\gamma^+}{\gamma^-}}}\alpha^\frac{\gamma^+}{\gamma^-}\\
&=\gamma^+\rho^{\gamma^+-\frac{\gamma^+}{\gamma^-}},
\end{align*}
  where we used
  \begin{align*}
  	\delta n=(1-\alpha)\rho_{-}(\rho,\delta n)=(1-\alpha)(\rho_{+}(\rho,\delta n))^{\frac{\gamma^+}{\gamma^-}}=(1-\alpha)(\dfrac{\rho}{\alpha})^{\frac{\gamma^+}{\gamma^-}}.
  \end{align*}

Finally, for the case   $n=\rho=0$, we have  $ \rho_+(0,0)= \rho_-(0,0)= P(0,0)=0$ and
\begin{align*}
		\partial_\rho \rho_+(0,0)&=	\lim _{\delta\rho \rightarrow 0^+}\frac{\rho_+(\delta\rho,0)}{\delta\rho}=1,\\
		\partial_\rho P(0,0)&=	\lim _{\delta\rho \rightarrow 0^+}\frac{P(\delta\rho,0)}{\delta\rho}=0,\\
		\partial_n \rho_+(0,0)&=	\lim _{\delta n \rightarrow 0^+}\frac{\rho_+(0,\delta n)}{\delta n}=\lim _{\delta n \rightarrow 0^+}\frac{\mathfrak{q}^{-1}(\delta n)}{\delta n}\\
		&=	\lim _{\delta n \rightarrow 0^+} (\delta n)^{\frac{\gamma^-}{\gamma+}-1}= \left\{\begin{array}{l}
		0, ~~\text{if}~~\frac{\gamma^-}{\gamma+}>1,\\
		1,~~\text{if}~~\frac{\gamma^-}{\gamma+}=1,
		\end{array}\right.\\
	\partial_n  P(0,0)&=\lim _{\delta n\rightarrow 0^+}\frac{P(0, \delta n)}{\delta n}= \lim _{\delta n\rightarrow 0^+}( \delta n)^{\gamma^--1}= \left\{\begin{array}{l}
		0, ~~\text{if}~~\gamma^->1,\\
		1,~~\text{if}~~\gamma^-=1,
	\end{array}\right.
\end{align*}
where we used $ \rho_+(\delta\rho,0)=\delta\rho, \rho_+(0,\delta n)= \mathfrak{q}^{-1}(\delta n)$.

In summary,  we have arrived at the conclusion.
 \begin{flushright}
      $\square$
      \end{flushright}

  \begin{lem} \label{le3}
  	If $\rho, n>0$,  then 
  	 $ n \mapsto \partial_{n} P(\rho, n)$   are Lipschitz continuous. Moreover,
  	  \begin{align}
  	 \left|\partial^2_{n} P(\rho, n)\right|\leq C\left(\rho^{\gamma^+-2\frac{\gamma^+}{\gamma^-}}+n^{\gamma^--2}\right),\label{P1}	
  	 \end{align}
   where  $C$ is a  positive constant.
  \end{lem}
  
  \noindent{\bf{Proof:}}
  
  We only need to prove \eqref{P1} and the analysis is similar with those of Lemma \ref{le1}.
  Direct calculation gives
  \begin{align}
  	\partial^2_{n} P(\rho, n)&=\mathfrak{P}''_{+}\left(\rho_{+}(\rho, n)\right)(\partial_{n} \rho_{+}(\rho, n))^2+\mathfrak{P}'_{+}\left(\rho_{+}(\rho, n)\right)\partial^2_{n} \rho_{+}(\rho, n)
  \end{align}
  and, with \eqref{rho+n}, we can get 
  \begin{align}
  	\mathfrak{P}''_{+}\left(\rho_{+}(\rho, n)\right)(\partial_{n} \rho_{+}(\rho, n))^2& \leq C(\rho_+)^{(\gamma^+-2)}(\rho_+)^{-2(\frac{\gamma^+}{\gamma^-}-1)}\nonumber\\
  	&\leq  C(\rho+\mathfrak{q}^{-1}(n))^{\gamma^+-2\frac{\gamma^+}{\gamma^-}}\nonumber\\
  	&\leq C\left(\rho^{\gamma^+-2\frac{\gamma^+}{\gamma^-}}+n^{\gamma^--2}\right), \label{Pn21}
  \end{align}
  where we have used \eqref{rho+1}.
  
  Similarly, one has 
  \begin{align} 
  	0<\left|\partial^2_{n} \rho_{+}(\rho, n)\right|&=\frac{\rho_{+}\partial_{n} \rho_{+}}{\left[\rho \mathfrak{q}\left(\rho_{+}\right)+\rho_{+} \mathfrak{q}^{\prime}\left(\rho_{+}\right)\left(\rho_{+}-\rho\right)\right]^2}
  	\left|2\rho\mathfrak{q}\left(\rho_{+}\right)-2\rho\rho_+\mathfrak{q}'\left(\rho_{+}\right)
  	-(\rho_+)^2\mathfrak{q}''\left(\rho_{+}\right)(\rho_+-\rho)\right|	\nonumber\\
  	&=\frac{(\rho_{+})^3}{\left[\rho \mathfrak{q}\left(\rho_{+}\right)+\rho_{+} \mathfrak{q}^{\prime}\left(\rho_{+}\right)\left(\rho_{+}-\rho\right)\right]^2}
  	\left|\frac{2\alpha(1-\frac{\gamma^+}{\gamma^-})-\frac{\gamma^+}{\gamma^-}(\frac{\gamma^+}{\gamma^-}-1)(1-\alpha)}{\alpha+\frac{\gamma^+}{\gamma^-}\left(1-\alpha\right)}\right|\nonumber\\
  	&\leq \bar{q}(\frac{\gamma^+}{\gamma^-}+2)\left|\frac{\gamma^+}{\gamma^-}-1\right|\frac{\left(\rho_{+}\right)^{3}}{\rho_{+}^4( \mathfrak{q}^{\prime}\left(\rho_{+}\right))^2}\frac{1}{[\frac{\gamma^-}{\gamma^+}\alpha+(1-\alpha)]^2}\nonumber\\
  	&\leq C\bar{q}\bar{q}^2_1(\frac{\gamma^+}{\gamma^-}+2)|\frac{\gamma^+}{\gamma^-}-1|\left(\rho_+\right)^{1-2\frac{\gamma^+}{\gamma^-}}, \label{rn2}
  \end{align}	
  where we used $\partial_{n} \rho_{+}(\rho, n)=\dfrac{\left(\rho_{+}\right)^{2}}{\rho \mathfrak{q}\left(\rho_{+}\right)+\rho_{+} \mathfrak{q}^{\prime}\left(\rho_{+}\right)\left(\rho_{+}-\rho\right)}$ and $\mathfrak{q}(\cdot)=(\cdot)^{\frac{\gamma^+}{\gamma^-}}$. Thus,
  \begin{align}
  	\left|\mathfrak{P}'_{+}\left(\rho_{+}(\rho, n)\right)\partial^2_{n} \rho_{+}(\rho, n)\right|& \leq C \left(\rho_+\right)^{\gamma^+-2\frac{\gamma^+}{\gamma^-}}\nonumber\\
  	&\leq C\left(\rho^{\gamma^+-2\frac{\gamma^+}{\gamma^-}}+n^{\gamma^--2}\right), \label{Pn22}
  \end{align}
which  combine with \eqref{Pn21} gives \eqref{P1}. 
  
  \begin{flushright}
  	$\square$
  \end{flushright}

We introduce function $s(t,x)$ as
\begin{equation} \label{s}
	s:=\left\{\begin{array}{l}
		\dfrac{\rho}{n} ~~\text {   if } n>0 \\
		
		0~~~ \text {   if } n=0,
	\end{array}\right.
\end{equation}
and when  $(\rho,n)\in 	\mathscr{O}_{c_0}$, one has $0\leq s \leq c_0$.  The following lemma is based on function $s$, and is crucial in process of passing the $\epsilon$, $\delta$ to $0^+$ in Section \ref{s4.} and \ref{s5.}.
\begin{lem} \label{le2}
	If $ (\rho,n)\in 	\mathscr{O}_{c_0} , 1\leq \gamma^+,\gamma^-$, then for any  $0\leq s \leq c_0$, $P(ns, n)$ is a non-decreasing of function $n$ and can be decomposed  into
\begin{align}
	P(ns, n)=\frac{\underline{q}_1}{2}n^{\gamma^-} +\pi(n, s), \label{Prs}
\end{align}
where $\underline{q}_1$ is a positive constant  defined as \eqref{q11} and  $[0, \infty) \ni n \mapsto \pi(n, s)$ is a non-decreasing  map.
\end{lem}

\noindent{\bf{Proof:}}

From \eqref{q}, \eqref{rho+0}, and \eqref{rho+}, we know
\begin{align}
	\rho_{-} \mathfrak{q}^{-1}\left(\rho_{-}\right)-\mathfrak{q}^{-1}\left(\rho_{-}\right) n-\rho\rho_{-}=0. \label{rho-}
\end{align}
We can see $\rho,\rho^+,\mathfrak{q}$ are symmetry to $n,\rho^-,\mathfrak{q}^{-1}$ and thus can follow the analysis in Lemma \ref{le1} to obtain that 
\begin{align}
\underline{q}_1\leq\partial_{n} \rho_{-}(\rho, n)\leq \overline{q}_1,& \label{rr1}\\
	0\leq\partial_{\rho} \rho_{-}(\rho, n) \leq C& \left(n^{1-\frac{\gamma^-}{\gamma^+}}+\rho^{\frac{\gamma^+}{\gamma-}-1}\right), \label{rho+2'}
\end{align}
for $\rho>0,n>0$, where $ \underline{q}_1,\overline{q}_1$ is defined as \eqref{q11}. If $ 1\leq \gamma^+,\gamma^- $, we can also  verify that \eqref{Pr}, \eqref{Pn} hold in the case $ \rho=0 $ or $ n=0 $.

We set
\begin{align}
	\pi(n, s)&:=P(ns, n)-\frac{\underline{q}_1}{2}n^{\gamma^-}\nonumber\\
	&=\mathfrak{P}_{-}\left(\rho_{-}(ns, n)\right)-\frac{\underline{q}_1}{2}n^{\gamma^-}.
\end{align}
If $n>0$, then $\rho>0$ due to $ n\leq c_0 \rho  $, and it follows \eqref{rr1}, \eqref{rho+2'} that
\begin{align}
	\frac{\partial \pi}{\partial n}&=\mathfrak{P}'_{-}\left(\rho_{-}(ns, n)\right)\left(s\frac{\partial \rho_{-}}{\partial \rho}(ns, n)+ \frac{\partial \rho_{-}}{\partial n}(ns, n)\right)-\frac{\underline{q}_1\gamma^{-}}{2}n^{\gamma^--1}\\
&\geq\gamma^{-}\left(\rho_{-}(ns, n)\right)^{\gamma^--1}\left(\underline{q}_1+0\right)-\frac{\underline{q}_1\gamma^{-}}{2}n^{\gamma^--1}\nonumber\\
	&\geq \gamma^{-}\underline{q}_1n^{\gamma^--1}-\frac{\underline{q}_1\gamma^{-}}{2}n^{\gamma^--1}\nonumber\\
	&\geq \frac{\underline{q}_1\gamma^{-}}{2}n^{\gamma^--1}\geq0, 
\end{align}
where we used $n=(1-\alpha)\rho_-\leq \rho_-$.   
   If $n=0$ and  $ 1\leq \gamma^+,\gamma^- $, it follows  \eqref{Pr}, \eqref{Pn} that
   \begin{align*}
   		\frac{\partial \pi}{\partial n}(0,s)&=s\partial_{\rho}P(0,0) +\partial_{n}P(0,0)-\frac{\underline{q}_1\gamma^-}{2}n^{\gamma^--1}|_{n=0}\geq 0.
   \end{align*}

In conclusion,  $	\dfrac{\partial \pi}{\partial n}(n, s)\geq 0$ for all $n\in[0, \infty)$, and $ 	P(ns, n) $ can be written into form of \eqref{Prs} and is non-decreasing, as $n^{\gamma^-}$ is increasing for $n\in [0,\infty)$.
\begin{flushright}
	$\square$
\end{flushright}

\begin{lem} \label{le4}
	Suppose $(\rho,n) \in 	\mathscr{O}_{c_0} \cap(\underline{r}, \infty)$, then
\begin{align}
	\left|\frac{\partial^{2} H_P}{\partial \rho^{2}}(\rho, n)\right|+\left|\frac{\partial^{2} H_P}{\partial \rho \partial n}(\rho, n)\right|+\left|\frac{\partial^{2} H_P}{\partial n^{2}}(\rho, n)\right| \leq C(\underline{r})\left(1+\rho^{A}\right), \label{H0}
\end{align}
where $\underline{r}>0$, $C(\underline{r})$ is a constant that depend on $\underline{r}$ and
\begin{align}
	A=\max\{0,\gamma^{+}-2,\gamma^{-}-2,\gamma^+-2\frac{\gamma^+}{\gamma^-},\gamma^+-\frac{\gamma^+}{\gamma^-}-1,\gamma^--\frac{\gamma^-}{\gamma^+}-1\}. \label{A}
\end{align}

\end{lem}

\noindent{\bf{Proof:}}

Recalling the definition of $H:=H_{P}(\rho, n)$ in \eqref{ener3} and  taking derivation on $H$ with respect to $\rho$, $n$, we have 
\begin{align}
&\frac{\partial H}{\partial\rho}=\int_{1}^{\rho} \frac{P\left(\zeta, \zeta \frac{n}{\rho}\right)}{\zeta^{2}} \mathrm{d} \zeta-\frac{n}{\rho} \int_{1}^{\rho} \frac{P'\left(\zeta, \zeta \frac{n}{\rho}\right)}{\zeta} \mathrm{d} \zeta+\frac{P(\rho, n)}{\rho},\\
&\frac{\partial H}{\partial n}=\int_{1}^{\rho}\frac{P'\left(\zeta, \zeta \frac{n}{\rho}\right)}{\zeta} \mathrm{d} \zeta,\\
&\frac{\partial^2 H}{\partial\rho\partial n}=-\frac{n}{\rho^2} \int_{1}^{\rho} P''\left(\zeta, \zeta \frac{n}{\rho}\right) \mathrm{d} \zeta+\frac{P'(\rho, n)}{\rho},\label{Hrn}\\
&\frac{\partial^2 H}{\partial n^2}=\frac{1}{\rho}\int_{1}^{\rho}P''\left(\zeta, \zeta \frac{n}{\rho}\right) \mathrm{d} \zeta,\label{Hnn}\\
&\frac{\partial^2 H}{\partial\rho^2}=-\frac{n}{\rho^2}P'(\rho, n)+\frac{\partial_{\rho} P(\rho,n)}{\rho}+\frac{n^2}{\rho^3} \int_{1}^{\rho} P''\left(\zeta, \zeta \frac{n}{\rho}\right)\mathrm{d}\zeta,\label{H2r}
\end{align}
where we denote $P'(\rho, n):=\partial_n P\left(\rho, n\right),~P''(\rho,n):=\partial^2_{n} P\left(\rho, n\right)$. By $ \rho\geq \underline{r} $, \eqref{P1} and \eqref{Hnn}, one can direct  calculate that
\begin{align}
	\begin{aligned}
	\left|\frac{\partial^2 H}{\partial n^2}\right|&\leq\frac{C}{\rho}\left|\int_{1}^{\rho}\left(\zeta^{\gamma^+-2\frac{\gamma^+}{\gamma^-}}+\zeta^{\gamma^--2}\frac{n^{\gamma^--2}}{\rho^{\gamma^--2}}\right)\mathrm{d} \zeta\right|\\
	&\leq \frac{C}{\rho}\left(\rho^{\gamma^+-2\frac{\gamma^+}{\gamma^-}+1}+1+\rho n^{\gamma^--2}+\frac{n^{\gamma^--2}}{\rho^{\gamma^--2}}\right)\\
	&\leq C\left(\rho^{\gamma^+-2\frac{\gamma^+}{\gamma^-}}+n^{\gamma^--2}\right)+\frac{C}{\underline{r}},
\end{aligned}\label{Hnn1}
\end{align}
which 
  combined with \eqref{Pn}, \eqref{Pr}, \eqref{P1},  \eqref{Hrn} and \eqref{H2r}  shows
\begin{align}
\left|\frac{\partial^2 H}{\partial\rho\partial n}\right|&\leq C\left(n\rho^{\gamma^+-2\frac{\gamma^+}{\gamma^-}-1}+\rho^{\gamma^+-\frac{\gamma^+}{\gamma^-}-1}+\frac{n^{\gamma^--1}}{\rho}\right)+\frac{C}{\underline{r}},\label{Hrn1}\\
		\left|\frac{\partial^2 H}{\partial \rho^2}\right|\leq C&\left(n^2\rho^{\gamma^+-2\frac{\gamma^+}{\gamma^-}-2}+n\rho^{\gamma^+-\frac{\gamma^+}{\gamma^-}-2}+\frac{n^{\gamma^-}}{\rho^2}+\rho^{\gamma^+-2}+\frac{n^{\gamma^--\frac{\gamma^-}{\gamma^+}}}{\rho}\right)+\frac{C}{\underline{r}}. \label{Hrr1} 
\end{align}
Combining \eqref{Hnn1}, \eqref{Hrn1}, \eqref{Hrr1} and using $\dfrac{\underline{r}}{c_0}<\dfrac{1}{c_0} n\leq \rho$, we can get \eqref{H0}.
\begin{flushright}
	$\square$
\end{flushright}

\subsection{Analysis tools}

~~~~In this subsection, we introduce some propositions and lemmas that are necessary for the analysis in later Sections.

The following lemma concerns the renormalized solutions introduced by DiPerna et al. \cite{diperna1989ordinary} and is presented in Vasseur et al. \cite{vasseur2019global}. When $N=1$, it is the corresponding result in Feireisl et al. \cite{feireisl2001existence}.

\begin{lem}\label{le2.1} \cite{vasseur2019global}~
Let $\beta: \mathbb{R}^{N} \rightarrow \mathbb{R}$ be a $C^{1}$ function with $|\nabla \beta(X)| \in L^{\infty}\left(\mathbb{R}^{N}\right),$ and $R \in\left(L^{2}\left(0, T ; L^{2}(\Omega)\right)\right)^{N}$
	$u \in L^{2}\left(0, T ; H_{0}^{1}(\Omega)\right)$ satisfy
\begin{align}
		\frac{\partial}{\partial_{t}} R+\operatorname{div}(\mathbf{u} \otimes R)=0,\left.\quad R\right|_{t=0}=R_{0}(x) \label{normal}
\end{align}
	in the distribution sense. Then we have
	$$
	(\beta(R))_{t}+\operatorname{div}(\beta(R)\mathbf{u} )+[\nabla \beta(R) \cdot R-\beta(R)] \operatorname{div}\mathbf{u} =0
	$$
	in the distribution sense. Moreover, if $R \in L^{\infty}\left(0, T ; L^{\gamma}(\Omega)\right)$ for $\gamma>1$, then
	$$
	R \in \left(C\left([0, T] ; L^{1}(\Omega)\right)\right)^N
	$$
	and so
	$$
	\int_{\Omega} \beta(R) d x(t)=\int_{\Omega} \beta\left(R_{0}\right) d x-\int_{0}^{t} \int_{\Omega}[\nabla \beta(R) \cdot R-\beta(R)] \operatorname{div} \mathbf{u}  d x d t.
	$$
\end{lem}
 \begin{flushright}
	$\lozenge$
\end{flushright}

The following proposition is standard in the analysis for obtaining weak solution for fluid models.
\begin{prop} \label{prop2.2} \cite{novotny2020weak}
 Let $P, G:(I \times \Omega) \times[0, \infty) \mapsto R$ be a couple of functions such that for almost all $(t, x) \in I \times \Omega, \rho \mapsto P(t, x, \rho)$ and $\rho \mapsto G(t, x, \rho)$ are both non
	decreasing and continuous on $[0, \infty) .$ Assume that $\rho_{n} \in L^{1}\left(Q_{T}\right)$ is a sequence such that
	$$
	\left.\begin{array}{c}
		P\left(\cdot, \rho_{n}(\cdot)\right) \rightarrow \overline{P(\cdot, \rho)} \\
		G\left(\cdot, \rho_{n}(\cdot)\right) \rightarrow \overline{G(\cdot, \rho)} \\
		P\left(\cdot, \rho_{n}(\cdot)\right) G\left(\cdot, \rho_{n}(\cdot)\right) \rightarrow \overline{P(\cdot, \rho) G(\cdot, \rho)}
	\end{array}\right\} \text { in } L^{1}(I \times \Omega).
	$$
	Then
	$$
	\overline{P(\cdot, \rho)} ~\overline{G(\cdot, \rho)} \leqq \overline{P(\cdot, \rho) G(\cdot, \rho)}
	$$
	a.e. in $I \times \Omega$.
\end{prop}
 \begin{flushright}
	$\lozenge$
\end{flushright}
\begin{prop} \label{prop2.3} \cite{vasseur2019global}
	Let $\varsigma_k \rightarrow 0$, as $k\rightarrow 0$. If $\rho_k\geq0$, $n_k\geq0$ are the solutions to
\begin{align}
	\partial_{t} \rho_k+\operatorname{div}(\rho_k \mathbf{u}_k) = \varsigma_k \Delta \rho_k \label{app50}\\
	\partial_{t} n_k+\operatorname{div}(n_k \mathbf{u}_k) =\varsigma_k \Delta n_k \label{app60}\\
	\left.\frac{\partial \rho_k}{\partial\nu}\right|_{\partial \Omega} =0,~ \left.\frac{\partial n_k}{\partial\nu}\right|_{\partial \Omega} =0 \label{app610}\\
	\rho_{k}(0, x) =\rho_{0}, ~n_{k} (0, x) =n_{0}. \label{app620}
\end{align}	
and
 $$
	\mathbf{u}_{k} \in L^{2}\left(I, W_{0}^{1,2}\left(\Omega \right)\right),\left(\rho_{k}, n_{k}\right) \in \mathcal{O}_{c_0} \cap\left(C\left(\bar{I} ; L^{1}(\Omega)\right) \cap L^{2}\left(Q_{T}\right)\right)^{2},
	$$
	$$
	\begin{array}{l}
		\sup _{k \in N}\left(\left\|\rho_{k}\right\|_{L^{\infty}\left(I ; L^{\gamma}(\Omega)\right)}+\left\|n_{k}\right\|_{L^{\infty}\left(I ; L^{\gamma}(\Omega)\right)}\right. \\
		\left.+\left\|\rho_{k}\right\|_{L^{2}\left(Q_{T}\right)}+\left\|\mathbf{u}_{k}\right\|_{L^{2}\left(I ; W^{1,2}(\Omega)\right)}\right)<\infty
	\end{array}
	$$
	and
	$$
	\int_{\Omega} n_{k}(0, x) s_{k}^{2}(0, x) \mathrm{d} x \rightarrow \int_{\Omega} n(0, x) s^{2}(0, x) \mathrm{d} x,
	$$
	where 	
	\begin{equation} 
		s_k:=\left\{\begin{array}{l}
			\dfrac{\rho_k}{n_k} ~~\text {   if } n>0 \\
			
			0~~~ \text {   if } n=0,
		\end{array}\right. 
	\quad
		s:=\left\{\begin{array}{l}
		\dfrac{\rho}{n} ~~\text {   if } n>0 \\
		
		0~~~ \text {   if } n=0,
	\end{array}\right.
	\end{equation}
and $$\frac{1}{c_{0}} n_{k} \leq \rho_{k}\leq c_{0} n, ~~\frac{1}{c_{0}} n \leq \rho_{k}\leq c_{0} n.$$ 
	
	Then, up to a subsequence, we have
\begin{align*}
	n_{k} &\rightarrow n, \quad \rho_{k} \rightarrow \rho \text { weakly in } L^{\infty}\left(0, T ; L^{2}(\Omega)\right), \\
& \mathbf{u}_{k} \rightarrow u \text { weakly in } L^{2}\left(0, T ; H_{0}^{1}(\Omega)\right)
\end{align*}
	and for any $p\geq1$,
	\begin{align}
	\int_{\Omega}\left(n_{k}\left|s_{k}-s\right|^{p}\right)(\tau, \cdot) \mathrm{d} x \rightarrow 0 \text { with any } 1 \leq p<\infty	 \label{appt}
	\end{align}
for all $\tau \in [0,T]$. 
\end{prop}

\begin{rem}
	Proposition \ref{prop2.3} is analogous to Theorem 2.2 in Vasseur et al. \cite{vasseur2019global} and Proposition 11 in Novotn\'y et al \cite{novotny2020weak}, and is the key to obtain strong convergence of $n_k$ and to handle the possible oscillation for the pressure law function.  We refer the proof to  \cite{vasseur2019global}.  When $\varsigma_k \equiv 0$, it can be applied in Sections \ref{s4.} and \ref{s5.} corresponding to the process $\epsilon\rightarrow 0^{+}$ and $\delta\rightarrow 0^{+}$.
	
	Corresponding to mono-fluid case in Feireisl et al. \cite{feireisl2001existence}, the best expected range of $\gamma^-$ would be $\gamma^->\frac{3}{2}$ and it is still an open problem for $bi\text{-}fluid$ to reach this range. The sufficient condition for the proof of  Proposition \ref{prop2.3} is that $ n_k$ satisfy the renormalized transport equation \eqref{normal}, and even with the help of estimate in Lemma \ref{le5.1} given by the Bogovskii operater, we only have $\gamma^-+\frac{2}{3}\gamma^--1\geq2$, i.e. $\gamma^-\geq \frac{9}{5}$. See more in \cite{novotny2020weak}.
\end{rem}

\section{Faedo-Galerkin approach} \label{s3.}

Motivated by the work of \cite{novotny2020weak} and \cite{hu2010global}, we propose  following approximation system
\begin{align} 
\partial_{t} \rho&+\operatorname{div}(\rho \mathbf{u}) = \epsilon \Delta \rho, \label{app1}\\
\partial_{t} n&+\operatorname{div}(n \mathbf{u}) =\epsilon \Delta n, \label{app2}\\
\partial_{t}((\rho+n) \mathbf{u})+\operatorname{div}((\rho+n) \mathbf{u} \otimes \mathbf{u})&+\nabla \Pi_{\delta}(\rho, n)=\mu \Delta \mathbf{u}+\epsilon\nabla \mathbf{u}\cdot\nabla (\rho+n)\nonumber\\
&+(\mu+\lambda) \nabla \operatorname{div} \mathbf{u}+(\nabla \times \mathbf{H}) \times \mathbf{H}, \label{app3}\\
\mathbf{H}_{t}-\nabla \times(\mathbf{u} \times \mathbf{H})=-&\nabla \times(\nu \nabla \times \mathbf{H}),~~~
\operatorname{div} \mathbf{H}=0 ,\label{app4}
\end{align}
on $(t,x) \in Q_T$, with initial boundary condition
\begin{align} 
	\left.(\rho, n,(\rho+n)\mathbf{u},\mathbf{H})\right|_{t=0}&=\left(\rho_{0, \delta}, n_{0, \delta}, \mathbf{m}_{0, \delta},\mathbf{H}_{0, \delta}\right) \text { on } \overline{\Omega} \label{app40}\\
\left(\frac{\partial \rho}{\partial\nu}\right.&, \left. \left.\frac{\partial n}{\partial\nu},\mathbf{u},\mathbf{H}\right)\right|_{\partial \Omega} =0, \label{app41}
\end{align}
where $\epsilon,\delta>0$, and 
\begin{equation}
\Pi_{\delta}(\rho, n):=P(\rho, n)+\delta\left(\rho^{B}+n^{B}+\frac{1}{2} \rho^{2} n^{B-2}+\frac{1}{2} n^{2} \rho^{B-2}\right),	\label{Pd}
\end{equation}
and $B$ is large enough such that $B\geq A+2$. Here, $A$ is defined in \eqref{A}.
The initial boundary condition satisfy $ \rho_{0, \delta},n_{0, \delta} \in C^{3}(\Omega), \mathbf{u}_{0, \delta},\mathbf{H}_{0, \delta} \in C_{0}^{3}(\Omega),$ $	\mathbf{m}_{0, \delta}=\left(\rho_{0, \delta}+n_{0, \delta}\right) \mathbf{u}_{0, \delta}$
and
\begin{equation} \label{app42}
	\begin{aligned}
			\frac{1}{c_{0}} \rho_{0, \delta}& \leq n_{0, \delta} \leq c_{0} \rho_{0, \delta} ,~~0<\delta \leq \rho_{0, \delta}, n_{0, \delta} \leq \delta^{-\frac{1}{2 B}},\left.\quad\left(\frac{\partial n_{0, \delta}}{\partial \nu}, \frac{\partial \rho_{0, \delta}}{\partial \nu}\right)\right|_{\partial \Omega}=0, \\
&	\lim _{\delta \rightarrow 0}\left(\left\|\rho_{0, \delta}-\rho_{0}\right\|_{L^{\gamma^+}(\Omega)}+\left\|n_{0, \delta}-n_{0}\right\|_{L^{\gamma^-}(\Omega)}+\left\|\mathbf{H}_{0, \delta}-\mathbf{H}_{0}\right\|_{L^2(\Omega)}\right)=0, \\
	&	\mathbf{u}_{0 \delta}=\frac{\varphi_{\delta}}{\sqrt{\rho_{0, \delta}+n_{0, \delta}}} \eta_{\delta} *\left(\frac{\mathbf{m}_{0}}{\sqrt{\rho_{0}+n_{0}}}\right) ,~~	\mathbf{m}_{0, \delta} \rightarrow \mathbf{m}_{0} \text { in } L^{1}(\Omega) \text { as } \delta \rightarrow 0, 
\end{aligned}
\end{equation}
where $ \eta_{\delta}$ is the standard mollifier.

To obtain the existence of solutions for \eqref{app1}-\eqref{app41}, we need to investigate  the properties of the magnetic field. Indeed, for system 
\begin{align} \label{H}
\left\{\begin{array}{l}
	\mathbf{H}_{t}-\nabla \times(\mathbf{u} \times \mathbf{H})=-\nabla \times(\nu \nabla \times \mathbf{H}), \\
	\operatorname{div} \mathbf{H}=0, \\
	\mathbf{H}(x, 0)=\mathbf{H}_{0},\left.\quad \mathbf{H}\right|_{\partial \Omega}=0,
\end{array}\right.
\end{align}
we have following lemma that is presented in \cite{hu2010global}.
\begin{lem} \label{appH} \cite{hu2010global} 
 Assume that $\mathbf{u} \in C\left(0, T; C_{0}^{2}(\bar{\Omega} )\right)$ is a given velocity field, then there exists at most one function
	$$
	\mathbf{H} \in L^{2}\left(0, T ; H_{0}^{1}(\Omega)\right) \cap L^{\infty}\left(0, T; L^{2}(\Omega)\right)
	$$
	which solves \eqref{H} in the weak sense on $\Omega \times(0, T),$ and satisfies boundary and initial conditions in the sense of traces.	
Moreover, the operator $\mathbf{u} \mapsto \mathbf{H}[\mathbf{u}]$ maps bounded sets in $C\left(0, T; C_{0}^{2}(\bar{\Omega} )\right)$ into bounded
	subsets of
	$$
	Y:=L^{2}\left(0, T ; H_{0}^{1}(\Omega)\right) \cap L^{\infty}\left(0, T ; L^{2}(\Omega)\right)
	$$
	and the mapping
	$$
	\mathbf{u} \in C\left(0, T ; C_{0}^{2}\left(\bar{\Omega} \right)\right) \mapsto \mathbf{H} \in Y
	$$
	is continuous on any bounded subsets of $C\left(0, T ; \bar{\Omega}\right)$.
\end{lem}

Next, consider a sequence of finite dimensional spaces
$$\mathbf{X}_{k}=\left[\operatorname{span}\left\{\phi_{j}\right\}_{j=1}^{k}\right]^{3}, \quad k \in\{1,2,3, \cdots\}$$
where $\left\{\phi_{i}\right\}_{i=1}^{\infty}$ is the set of the eigenfunctions of the Laplacian:
$$
\left\{\begin{array}{l}
	-\Delta \phi_{i}=\lambda_{i} \phi_{i} \quad \text { on } \Omega \\
	\left.\phi_{i}\right|_{\partial \Omega}=0.
\end{array}\right.
$$
For any given $\epsilon, \delta>0,$ we need to find the approximate solution $\mathbf{u}_{k} \in C\left(0, T ; \mathbf{X}_{k}\right)$  given by:
\begin{align}
	&\int_{\Omega}(\rho_{k}+n_{k}) \mathbf{u}_{k}\cdot{\bm\phi} d x-\int_{\Omega} \mathbf{m}_{0, \delta}\cdot{\bm\phi}d x=	\int_{0}^{T} \int_{\Omega}\left((\rho_{k}+n_{k})(\mathbf{u}_{k} \otimes \mathbf{u}_{k}): \nabla{\bm\phi}+\Pi_{\delta}(\rho_{k}, n_{k}) \operatorname{div} {\bm\phi}\right) \mathrm{d} x \mathrm{d} t \nonumber\\
	&-\int_{0}^{T} \int_{\Omega}\left(\mu \nabla \mathbf{u}_{k}: \nabla{\bm\phi}+(\mu+\lambda) \operatorname{div} \mathbf{u}_{k} \operatorname{div} {\bm\phi}+\epsilon\nabla(\rho_{k}+n_{k}) \cdot \nabla \mathbf{u}_{k} \cdot {\bm\phi}-(\nabla \times \mathbf{H}_{k}) \times \mathbf{H}_{k} )\cdot{\bm\phi}\right) \mathrm{d} x \mathrm{d} t \label{app0}
\end{align}
for ${\bm\phi} \in C\left(0, T; \mathbf{X}_{k}\right)$, where $\mathbf{H}_{k}=\mathbf{H}_{k}\left(\mathbf{u}_{k}\right)$ satisfy \eqref{H}, and $ \rho_{k}=\rho_{k}\left(\mathbf{u}_{k}\right),~n_{k}=n_{k}\left(\mathbf{u}_{k}\right)$  satisfy
\begin{align}
\partial_{t} \rho_k+\operatorname{div}(\rho_k \mathbf{u}_k) = \epsilon \Delta \rho_k \label{app5}\\
\partial_{t} n_k+\operatorname{div}(n_k \mathbf{u}_k) =\epsilon \Delta n_k \label{app6}\\
\left.\frac{\partial \rho_k}{\partial\nu}\right|_{\partial \Omega} =0,~ \left.\frac{\partial n_k}{\partial\nu}\right|_{\partial \Omega} =0 \label{app61}\\
\rho_{k}(0, x) =\rho_{0}, ~n_{k} (0, x) =n_{0}. \label{app62}
\end{align}

With Lemma \ref{appH} and Lemma 2.2 of \cite{feireisl2001existence}, the problem \eqref{app0}, \eqref{H}, \eqref{app5}-\eqref{app62} can be solved locally
in time by means of the Schauder fixed-point argument, see Section
7 of Feireisl \cite{feireisl2004dynamics}.

Next, we need to derive the energy inequality to give the uniform estimates to extend $T_k$ to $T$. We can build the relationship between $n_k$ and $\rho_k$. Indeed, we have following lemma 
\begin{lem} \label{le3.2} \cite{vasseur2019global}
	If $\left(\rho_{k}, n_{k}, \mathbf{u}_{k}\right)$ is a solution to \eqref{app5} and \eqref{app6} with the initial data satisfying
	$$
\frac{1}{c_{0}} n_{0} \leq \rho_{0} \leq c_{0} n_{0}, \text{   on } \Omega
$$
for a given positive constant $c_0\geq 1$, then the following inequality holds
		\begin{equation}
		\frac{1}{c_{0}} n_{k}(x, t) \leq \rho_{k}(x, t) \leq c_{0} n_{k}(x, t) \text{   on } Q_T. \label{appest0}
	\end{equation}
\end{lem}
This could be done immediately by the maximum principle of parabolic equation. We refer the proof to \cite{vasseur2019global}.

Now, we set $ \delta>1 $, then by Lemma \ref{le3.2} we have $ \rho,n>1 $. Multiplying \eqref{app5}, \eqref{app6} by $\rho_k$ and $n_k$ respectively, integrating them on $\Omega$, one has
\begin{align}
\frac{1}{2} \frac{\mathrm{d}}{\mathrm{d} t}\int_{\Omega}\left(\rho^{2}_k+ n^{2}_k\right)dx+\epsilon\int_{\Omega}&\left(|\nabla \rho_k|^{2}+ |\nabla n_k|^{2}\right)dx
=-\frac{1}{2} \int_{\Omega}\left(\rho^{2}_k+n^{2}_k\right) \operatorname{div} \mathbf{u}_{k} d x \nonumber\\
&\leq C(\Sigma(\delta),\mu)\int_{\Omega}\left(\rho^{4}_k+n^{4}_k\right)dx +\frac{\mu}{4\Sigma(\delta)}\int_{\Omega}|\nabla \mathbf{u}_{k}|^2 d x,  \label{app7}
\end{align}
where $\Sigma(\delta)$ is a constant to be determined later. Take ${\bm\phi}=\mathbf{u}_k$ in \eqref{app0} and combine with \eqref{app5} and \eqref{app6} to get
\begin{align}
\frac{1}{2} \frac{\mathrm{d}}{\mathrm{d} t}\int_{\Omega} \left( (\rho_k+ n_k)|\mathbf{u}_k|^2 +|\mathbf{H}_k|^2\right)dx &+\int_{\Omega}\left(\mu|\nabla \mathbf{u}_k|^{2}+(\lambda+\mu)(\operatorname{div} \mathbf{u}_k)^{2}\right.\nonumber\\
&\left.+\nu|\nabla \times \mathbf{H}_k|^{2}\right) \mathrm{d} x= \int_{\Omega} \Pi_{\delta}(\rho_k, Z_k))\operatorname{div} \mathbf{u}_kdx.  \label{app8}
 \end{align}
Denote
 \begin{align}
\mathcal{H}_{\delta}(\rho, n):=H_{P}(\rho, n)+h_{\delta}(\rho, n), \label{H1}
\end{align}
where $H_{P}(\rho, n)$ is defined as \eqref{ener3} and
\begin{align}
h_{\delta}(\rho, n):=\frac{\delta}{B-1}\left(\rho^{B}+n^{B}+\frac{1}{2} \rho^{2} n^{B-2}+\frac{1}{2} n^{2} \rho^{B-2}\right). \label{h}	
\end{align}
Then, direct calculation gives
\begin{equation}
\begin{aligned}
\int_{\Omega}\partial_t\mathcal{H}_{\delta}(\rho_k, n_k)dx&=\int_{\Omega}\left(\partial_{\rho}\mathcal{H}_{\delta}(\rho_k, n_k) \partial_t\rho_k+\partial_{n}\mathcal{H}_{\delta}(\rho_k, n_k) \partial_tn_k\right)dx,\\
\int_{\Omega} \left(\epsilon \Delta \rho_k \partial_{\rho}\mathcal{H}_{\delta}\right.&\left.+\epsilon \Delta n_k \partial_{n}\mathcal{H}_{\delta}\right)dx=
-\epsilon\int_{\Omega}\left(\frac{\partial^{2} \mathcal{H}_{\delta}}{\partial \rho^{2}}\left|\nabla \rho_k\right|^{2}\right.\\
&\left.+2 \frac{\partial^{2} \mathcal{H}_{\delta}}{\partial \rho \partial n} \nabla \rho_k \cdot \nabla n_k
+\frac{\partial^{2} \mathcal{H}_{\delta}}{\partial n^{2}}\left|\nabla n_k\right|^{2}\right) \mathrm{d} x,\\
\int_{\Omega} \left[\operatorname{div}(\rho_k \mathbf{u}_k) \partial_\rho\mathcal{H}_{\delta}\right.&+\left.\operatorname{div}(n_k\mathbf{u}_k) \partial_n\mathcal{H}_{\delta}\right]dx =\int_{\Omega} \left[ n_k\partial_n\mathcal{H}_{\delta}+ \rho_k\partial_\rho\mathcal{H}_{\delta}-\mathcal{H}_{\delta}\right]\operatorname{div} \mathbf{u}_kdx\\
&=\int_{\Omega} \left( n_k\partial_n+ \rho_k\partial_\rho-1\right)
\left(H_{P_{\delta}}+h_{\delta}\right)\operatorname{div} \mathbf{u}_kdx\\
&=\int_{\Omega} \Pi_{\delta}(\rho_k, n_k))\operatorname{div} \mathbf{u}_kdx,
\end{aligned}   \label{app90}
\end{equation}
where $\mathcal{H}_{\delta}=\mathcal{H}_{\delta}(\rho_k, n_k)$.
 Multiplying \eqref{app5}, \eqref{app6} with $\partial_{\rho}\mathcal{H}_{\delta}$ and $\partial_{n}\mathcal{H}_{\delta}$ respectively, integrating them on $\Omega$, combining them with \eqref{app90}, and we finally get
\begin{align}
\frac{\mathrm{d}}{\mathrm{d} t}\int_{\Omega}\mathcal{H}_{\delta}(\rho_k, n_k)dx&+\int_{\Omega} \Pi_{\delta}(\rho_k, n_k)\operatorname{div} \mathbf{u}_kdx=
-\epsilon\int_{\Omega}\left(\frac{\partial^{2} \mathcal{H}_{\delta}}{\partial \rho^{2}}\left|\nabla \rho_k\right|^{2}\right.\nonumber\\
&\left.+2 \frac{\partial^{2} \mathcal{H}_{\delta}}{\partial \rho \partial n} \nabla \rho_k\cdot \nabla n_k
+\frac{\partial^{2} \mathcal{H}_{\delta}}{\partial n^{2}}\left|\nabla n_k\right|^{2}\right) \mathrm{d} x.  \label{app9}
\end{align}

Multiplying \eqref{app7} with $\Sigma=\Sigma(\delta)$, combining it with \eqref{app8} and \eqref{app9}, we can obtain
\begin{align}
\frac{1}{2}& \frac{\mathrm{d}}{\mathrm{d} t}\int_{\Omega}\left(\Sigma(\rho^{2}_k+ n^{2}_k)+(\rho_k+ n_k)|\mathbf{u}_k|^2 +|\mathbf{H}_k|^2+2\mathcal{H}_{\delta}(\rho_k, n_k)\right)dx+\nonumber\\
&\epsilon\Sigma\int_{\Omega}\left(|\nabla \rho_k|^{2}+ |\nabla n_k|^{2}\right)dx+\int_{\Omega}\left(\frac{3\mu}{4}|\nabla \mathbf{u}_k|^{2}+(\lambda+\mu)(\operatorname{div} \mathbf{u}_k)^{2}+\nu|\nabla \times \mathbf{H}_k|^{2}\right) \mathrm{d} x\nonumber\\
&+\epsilon\int_{\Omega}\left(\frac{\partial^{2} h_{\delta}}{\partial \rho^{2}}\left|\nabla \rho_k\right|^{2}+2 \frac{\partial^{2} h_{\delta}}{\partial \rho \partial n} \nabla \rho_k\cdot \nabla n_k
+\frac{\partial^{2} h_{\delta}}{\partial n^{2}}\left|\nabla n_k\right|^{2}\right) \mathrm{d} x
\leq \nonumber\\
&-\epsilon\int_{\Omega}\left(\frac{\partial^{2} H_{P}}{\partial \rho^{2}}\left|\nabla \rho_k\right|^{2}+2 \frac{\partial^{2} H_{P}}{\partial \rho \partial n} \nabla \rho_k\cdot \nabla n_k
+\frac{\partial^{2} H_{P}}{\partial n^{2}}\left|\nabla n_k\right|^{2}\right) \mathrm{d} x \nonumber\\
&\quad\quad\quad+C(\Sigma,\mu)\Sigma \int_{\Omega}\left(\rho^{4}_k+n^{4}_k\right)dx. \label{app10}
\end{align}

If $\rho_k\geq 1,n_k\geq 1$, computing on \eqref{h}  and utilizing \eqref{appest0}, $ n>1 $ shows that
$$\left|\frac{\partial^{2} h_{\delta}(\rho_{k},n_k)}{\partial \rho^{2}}\right|,\left|\frac{\partial^{2} h_{\delta}(\rho_{k},n_k)}{\partial \rho \partial n}\right|,\left|\frac{\partial^{2} h_{\delta}(\rho_{k},n_k)}{\partial n^{2}} \right| \leq C\frac{\delta}{B-1} \rho^{B-2}_k,$$
which combined with Lemma \ref{le4}  and \eqref{appest0} leads to
\begin{align}
\epsilon\int_{\Omega}1_{\left\{\rho_{k}+n_{k} \geq K\right\}}&\left(\left|\frac{\partial^{2} H_{P}}{\partial \rho^{2}}\right|\left|\nabla \rho_k\right|^{2}+2\left| \frac{\partial^{2}H_{P}}{\partial \rho \partial n} \right|\nabla \rho_k\cdot \nabla n_k
+\left|\frac{\partial^{2} H_{P}}{\partial n^{2}}\right|\left|\nabla n_k\right|^{2}\right) \mathrm{d} x\nonumber\\
&\leq \epsilon\int_{\Omega}1_{\left\{\rho_{k}+n_{k} \geq K\right\}}\left(1+\rho_k^{A}\right)\left(\left|\nabla \rho_k\right|^{2}+2  \nabla \rho_k\cdot \nabla n_k
+\left|\nabla n_k\right|^{2}\right) \mathrm{d} x\nonumber\\
&\leq \epsilon\int_{\Omega}1_{\left\{\rho_{k}+n_{k} \geq K\right\}}\left(\rho^{B-2}_k+n^{B-2}_k\right)\left(\left|\nabla \rho_k\right|^{2}+2  \nabla \rho_k\cdot \nabla n_k
+\left|\nabla n_k\right|^{2}\right) \mathrm{d} x
\nonumber\\
&\leq \frac{C_2\epsilon(B-1)}{\delta}\int_{\Omega}\left(\frac{\partial^{2} h_{\delta}}{\partial \rho^{2}}\left|\nabla \rho_k\right|^{2}+2 \frac{\partial^{2} h_{\delta}}{\partial \rho \partial n} \nabla \rho_k\cdot \nabla n_k
+\frac{\partial^{2} h_{\delta}}{\partial n^{2}}\left|\nabla n_k\right|^{2}\right) \mathrm{d} x, \label{app11}
\end{align}
and
\begin{align}
\epsilon\int_{\Omega}1_{\left\{\rho_{k}+n_{k} < K\right\}}&\left(\left|\frac{\partial^{2} H_{P}}{\partial \rho^{2}}\right|\left|\nabla \rho_k\right|^{2}+2 \left|\frac{\partial^{2} H_{}}{\partial \rho \partial n}\right| \nabla \rho_k\cdot \nabla n_k
+\left|\frac{\partial^{2}H_{P}}{\partial n^{2}}\right|\left|\nabla n_k\right|^{2}\right) \mathrm{d} x\nonumber\\
&\leq \epsilon C(\delta)\int_{\Omega}1_{\left\{\rho_{k}+n_{k} < K\right\}}\left(n_k^A+\rho_k^{A}\right)\left(\left|\nabla \rho_k\right|^{2}+2  \nabla \rho_k\cdot \nabla n_k
+\left|\nabla n_k\right|^{2}\right) \mathrm{d} x\nonumber\\
&\leq2\epsilon  C(\delta)K^A\int_{\Omega}1_{\left\{\rho_{k}+n_{k} < K\right\}}\left(\left|\nabla \rho_k\right|^{2}+2  \nabla \rho_k\cdot \nabla n_k
+\left|\nabla n_k\right|^{2}\right) \mathrm{d} x\nonumber\\
&\leq\frac{ \epsilon\Sigma(\delta)}{2}\int_{\Omega}\left(|\nabla \rho_k|^{2}+ |\nabla n_k|^{2}\right)dx  \label{app12}
\end{align}
for $B\geq A+2,~K\geq 1$ and 
\begin{align}
	\Sigma(\delta)\geq4C(\delta)K^A, \label{sigma}
\end{align}
where $C(\delta)$ is a constant defined in Lemma \eqref{le4}.
We find that  the last line of \eqref{app12} is absorbed by second line of \eqref{app10} and the last line of \eqref{app11} is absorbed by the third line of  \eqref{app10}, if $\delta\geq \frac{(B-1)}{C_2}$.

On the other hand, for $\rho_k, n_k$ satisfying
\begin{align}
\int_{\Omega}\left(\frac{\delta}{B-1}(\rho_k^{B}+n_k^{B})\right) \mathrm{d} x
&\leq C(\Sigma,\mu)\Sigma \int_0^{T}\int_{\Omega}\left(\rho^{4}_k+n^{4}_k\right)dx\nonumber\\
&\leq C(\Sigma,\mu)\Sigma \int_0^{T} \int_{\Omega}\left(\rho^{B}_k+n^{B}_k\right)dx,
\end{align}
we can use Gronwall lemma  to obtain
\begin{align}
\int_{\Omega}\left(\rho_k^{4}+n_k^{4}\right) \mathrm{d} \leq \int_{\Omega}\left(\rho_k^{B}+n_k^{B}\right) \mathrm{d} x
\leq C(\Sigma,\mu)\Sigma T \exp^{C(\Sigma,\mu)\Sigma  T}. \label{app13}
\end{align}
Therefore, integrating \eqref{app10} on $[0,T]$, applying Gronwall inequality and combining it \eqref{app11}-\eqref{app13}, we have
\begin{align}
& \int_{\Omega}\left(\Sigma(\delta)(\rho^{2}_k+ n^{2}_k)+(\rho_k+ n_k)|\mathbf{u}_k|^2 +|\mathbf{H}_k|^2+2\mathcal{H}_{\delta}(\rho_k, n_k)\right)\mathrm{d}x+\nonumber\\
&\int_0^{T}\int_{\Omega}\left(\mu|\nabla \mathbf{u}_k|^{2}+(\lambda+\mu)(\operatorname{div} \mathbf{u}_k)^{2}+\nu|\nabla \times \mathbf{H}_k|^{2}\right) \mathrm{d} x\mathrm{d}t+\epsilon\int_0^T\int_{\Omega}\left(|\nabla \rho_k|^{2}+ |\nabla n_k|^{2}\right)dx\mathrm{d}t+\nonumber\\
&+\epsilon\int_0^T\int_{\Omega}\left(\left|\nabla \rho_k\right|^{2}
+\left|\nabla n_k\right|^{2}\right) \left(\rho_k^{B-2}+n_k^{B-2}\right)\mathrm{d} x\mathrm{d}t
\leq E_\delta(0)+ C(\Sigma,\mu)\Sigma T \exp^{C(\Sigma,\mu)\Sigma  T},  \label{appest1}
\end{align}
where 
\begin{align} \label{Ed0}
E_\delta(t):=\int_{\Omega}\left(\Sigma(\delta)(\rho^{2}_k+ n^{2}_k)+(\rho_k+ n_k)|\mathbf{u}_k|^2 +|\mathbf{H}_k|^2+2\mathcal{H}_{\delta}(\rho_k, n_k)\right)\mathrm{d}x.	
\end{align}
We have obtained a solution  $(\rho_k,n_k,\mathbf{u}_k,\mathbf{H}_k)$ globally in time for system \eqref{app0}, \eqref{H}, \eqref{app5}-\eqref{app62} with the use of \eqref{appest1}, which gives bounds 
\begin{equation}
\begin{aligned}
&0<\frac{1}{c_{k}} \leq \rho_{k}(x, t), n_{k}(x, t) \leq c_{k} \text { for a.e. }(x, t) \in \Omega \times(0, T), \\
&\sup _{t \in[0, T]}\left\|\rho_{k}(t)\right\|_{L^{B}(\Omega)} \leq C ,~~\sup _{t \in[0, T]}\left\|n_{k}(t)\right\|_{L^{B}(\Omega)} \leq C, \\
&\sup_{t \in[0, T]}\left\|\mathcal{H}_{\delta}(\rho_k, n_k)\right\|_{L^1(\Omega)} \leq C,
~~\int_{0}^{T}\left\|\nabla \times \mathbf{H}_k (t)\right\|_{L^{2}(\Omega)}\mathrm{d}t\leq C, \\
&\sup _{t \in[0, T]}\left\|\mathbf{H}_k(t)\right\|_{L^{2}(\Omega)} \leq C , ~~ \sup _{t \in[0, T]}
 \| \sqrt{(\rho_{k}+n_{k}) }\mathbf{u}_{k}(t) \|_{L^{2}(\Omega)}^{2} \leq C,\\
 &\epsilon \int_{0}^{T}\left(\left\|(\rho_k^{\frac{B}{2}-1}\nabla \rho_{k}(t)\right\|_{L^{2}(\Omega)}^{2}+\left\|n_k^{\frac{B}{2}-1}\nabla n_{k}(t)\right\|_{L^{2}(\Omega)}^{2}\right)\mathrm{d}t \leq C,
\end{aligned} \label{est1}
\end{equation}
where $C=C\left(T,\delta,\rho_0,n_0,\mathbf{m}_0,\mathbf{H}_0\right)$.

With these estimates and Aubin-Lions lemma, we can let the $k$  go to infinity to obtain a global weak solution for the approximation system \eqref{app1}-\eqref{app42}. Indeed,  for the pressure term 
$$
\Pi_{\delta}(\rho_k, n_k)=P(\rho_k, n_k)+\delta\left(\rho_k^{B}+n_k^{B}+\frac{1}{2} \rho_k^{2} n_k^{B-2}+\frac{1}{2} n_k^{2} \rho_k^{B-2}\right).
$$
in \eqref{app0}, \eqref{est1} gives
$$ \left\|\rho^{\frac{B}{2}}_{k}(t)\right\|_{L^{2}\left(0, T ; H^{1}(\Omega)\right)}\leq C,~~\sup _{t \in[0, T]}\left\|\rho_{k}(t)\right\|_{L^{B}(\Omega)} \leq C.$$
It follows interpolation inequality and Sobolev embedding inequality that
$$\begin{aligned}
	\int_{0}^{T}\left\|\rho_{k}^{B}(t)\right\|_{L^{2}}^{4 / 3} d t & \leq C \int_{0}^{T}\left\|\rho_{k}^{B}\right\|_{L^{1}}^{1 / 3}\left\|\rho_{k}^{B}\right\|_{L^{3}} d t \\
	& \leq C \sup _{t \in[0, T]}\left\|\rho_{k}^{B}(t)\right\|_{L^{1}}^{1 / 3} \int_{0}^{T}\left\|\rho_{k}^{B}\right\|_{L^{3}} d t\\
	& \leq C \sup _{t \in[0, T]}\left\|\rho_{k}^{B}(t)\right\|_{L^{1}}^{1 / 3} \int_{0}^{T}\left\|\nabla\rho_{k}^{\frac{B}{2}}\right\|_{L^{2}} d t\leq C.
\end{aligned}$$
If $B$ large enough, we have 
\begin{equation}
	\|\rho_{k}\|_{L^{B+1}(Q_T)}\leq C, \label{est11}
\end{equation}
 and the same analysis gives \begin{equation}
 	\|n_{k}\|_{L^{B+1}(Q_T)}\leq C. \label{est12}
 \end{equation}
 Using the estimates in \eqref{est1} and applying Aubin-Lions lemma (see  \cite{novotny2004introduction}) on \eqref{app5} and \eqref{app6}, we get $$n_k\rightarrow n ,~~~\rho_k\rightarrow \rho \text{ in } L^2(Q_T).$$
 This with \eqref{est11}, \eqref{est12} and interpolation inequality give
 \begin{equation}
 	n_k\rightarrow n ,~~~\rho_k\rightarrow \rho \text{ in } L^{B}(Q_T).
 \end{equation}
 Thus, 
 \begin{align}
 \|	\rho_k^2 n_k^{B-2}&-\rho^2n^{B-2}\|_{L^1(Q_T)}\leq \|	(\rho_k^2 -\rho^2)n_k^{B-2}\|_{L^1(Q_T)}+ \|	\rho^2 (n_k^{B-2}-n^{B-2})\|_{L^1(Q_T)}\nonumber\\
 &\leq C \|	\rho_k -\rho\|_{L^{B}(Q_T)}\|\rho_k +\rho\|_{L^{B}(Q_T)} \|	n_k^{B-2}\|_{L^{\frac{B}{B-2}}(Q_T)}\nonumber\\
 &+\|\rho^2\|_{L^\frac{B}{2}}\left(\|n_k^{B-3}\|_{L^{\frac{B}{B-3}}(Q_T)} +\|n^{B-3}\|_{L^{\frac{B}{B-3}}(Q_T)} \right)\|n_k -n\|_{L^{B}(Q_T)}\nonumber\\
 &\rightarrow 0,~\text{as }k \rightarrow \infty,
 \end{align}
 where we have used
 \begin{align*}
  n_k^{B-2}-n^{B-2}&=\int^1_0(B-2)\left[\sigma n_k+(1-\sigma)n\right]^{B-3}\left(n_k-n\right)d\sigma\\
  &\leq	C\left(|n_k|^{B-3}+|n|^{B-3}\right)|n_k-n|.
 \end{align*}
In the same manner,
$$ n_k^{2} \rho_k^{B-2}\rightarrow n^{2} \rho^{B-2}\text{ in } L^1(Q_T).$$
For $0<\theta<1$, it follows mean value theorem that
\begin{align}
	\|P(\rho_{k}&,n_k)-P(\rho,n)\|_{L^1(Q_T)}\nonumber \\
	&
	=\|\partial_\rho P\left(\rho+\theta(\rho_{k}-\rho),n_k\right)(\rho_k-\rho)+\partial_n P\left(\rho,n+\theta(n_k-n)\right)(n_k-n)\|_{L^1(Q_T)}\nonumber\\
	&\leq C \left(\|\partial_\rho P\left(\rho+\theta(\rho_{k}-\rho),n_k\right)\|_{L^{\frac{B}{B-1}}(Q_T)} \|	\rho_k -\rho\|_{L^{B}(Q_T)}\right.\nonumber\\
	&
	+\left.\|\partial_nP\left(\rho,n+\theta(n_k-n)\right)\|_{L^{\frac{B}{B-1}}(Q_T)} \|	n_k -n\|_{L^{B}(Q_T)}\right)\rightarrow 0,~\text{as }k \rightarrow \infty,
\end{align}
where
\begin{align*}
\|\partial_\rho P&\left(\rho+\theta(\rho_{k}-\rho),n_k\right)\|_{L^{\frac{B}{B-1}}(Q_T)}\\
&\leq C \left\|(\theta\rho_k+(1-\theta)\rho)^{\gamma^+-1}+(\theta\rho_k+(1-\theta)\rho)^{\gamma^--\frac{\gamma^-}{\gamma^+}}\right\|_{L^{\frac{B}{B-1}}(Q_T)}	\leq C,\\
\|\partial_nP&\left(\rho,n+\theta(n_k-n)\right)\|_{L^{\frac{B}{B-1}}(Q_T)}\\
&\leq C \left\|(\theta n_k+(1-\theta)n)^{\gamma^+-\frac{\gamma^+}{\gamma^-}}+(\theta n_k+(1-\theta)n)^{\gamma^--1}\right\|_{L^{\frac{B}{B-1}}(Q_T)}	\leq C
\end{align*} 
 due to \eqref{Pr}, \eqref{P1}, \eqref{est11}, \eqref{est12} and $B\geq A$.

For the magnetic field,  \eqref{est1} gives $\mathbf{H}_{k}\in L^{\infty}\left([0, T] ; L^{2}(\Omega)\right) \cap L^{2}\left([0, T] ; H_{0}^{1}(\Omega)\right)$, which combined with \eqref{app4} and the Aubin-Lions lemma yields strong convergence of 
$$\mathbf{H}_{k} \rightarrow \mathbf{H}\in L^{\infty}\left([0, T] ; L^{2}(\Omega)\right) \cap L^{2}\left([0, T] ; H_{0}^{1}(\Omega)\right).$$
Thus, we have
$$
\begin{aligned}
\left(\nabla \times \mathbf{H}_{k}\right) \times \mathbf{H}_{k} \rightarrow(\nabla \times \mathbf{H}) \times \mathbf{H} \text { in } \mathcal{D}^{\prime}(\Omega \times(0, T)),\\
\nabla \times\left(\mathbf{u}_{k} \times \mathbf{H}_{k}\right) \rightarrow \nabla \times(\mathbf{u} \times \mathbf{H}) \text { in } \mathcal{D}^{\prime}(\Omega \times(0, T))
\end{aligned}
$$
for
$
\mathbf{u}_{k} \rightharpoonup \mathbf{u} \text { in } L^{2}\left([0, T] ; H_{0}^{1}(\Omega)\right).
$

For the dispose of other terms in \eqref{app1}-\eqref{app4}, we can follow the same path of \cite{feireisl2001existence}, as there is no essential difference. We remark that the $C$ above equal to $C\left(T, \delta,\rho_0,n_0,\mathbf{m}_0,\mathbf{H}_0\right)$ and is independent of $k$ and can let $k$ in \eqref{appest0} and \eqref{appest1} go to infinity. In consequence, we have obtained following proposition

\begin{prop} \label{prop2}
	If $B\geq A+2$,  $\delta\geq \frac{(B-1)}{C_2}$, where $ C_2 $ is determined in \eqref{app11},  then there exists a global weak solution  $(\rho,n,\mathbf{u},\mathbf{H})$ to \eqref{app1}-\eqref{app42} for any given $T>0$ and hold following estimates 
\begin{align}
	&\frac{1}{c_{0}} n(x, t) \leq \rho(x, t) \leq c_{0} n(x, t), \label{appest2}\\
 \int_{\Omega}\left(\rho^{2}+ n^{2}+\right.&\left.(\rho+ n)|\mathbf{u}|^2 +|\mathbf{H}|^2+2\mathcal{H}_{\delta}(\rho, n)\right)\mathrm{d}x+\nonumber\\
\int_0^{T}\int_{\Omega}\left(\mu|\nabla \mathbf{u}|^{2}+(\lambda+\mu)\right.&\left.(\operatorname{div} \mathbf{u})^{2}+\nu|\nabla \times \mathbf{H}|^{2}\right) \mathrm{d} x\mathrm{d}t+\epsilon\int_0^T\int_{\Omega}\left(|\nabla \rho|^{2}+ |\nabla n|^{2}\right)dx\mathrm{d}t+\nonumber\\
+\epsilon\int_0^T\int_{\Omega}\left(\left|\nabla \rho\right|^{2}
+\left|\nabla n\right|^{2}\right) &\left(\rho^{B-2}+n^{B-2}\right)\mathrm{d} x\mathrm{d}t
\leq E_\delta(0)+ C(\Sigma,\mu)\Sigma T \exp^{C(\Sigma,\mu)\Sigma  T}, \label{appest3}
\end{align}
where $E_\delta(t)$, $\Sigma=\Sigma(\delta)$ are defined in \eqref{Ed0}, \eqref{sigma} respectively, and $ C(\Sigma,\mu) $ is independent of $ \epsilon $.
\end{prop}

\section{The vanishing viscosity limit $\epsilon\rightarrow 0^{+}$} \label{s4.}

\quad In this section, we aim at passing to the limit of $ (\rho_\epsilon, n_\epsilon,\mathbf{u}_{\epsilon},\mathbf{H}_{\epsilon})$ as $\epsilon$
goes to zero. The method is similar with those in Feireisl et al. \cite{feireisl2001existence} and Novotn{\'y} et al. \cite{novotny2020weak}.  We need to build a better estimate for $\rho_\epsilon, n_\epsilon$, as those in \eqref{appest3} is not enough to obtain the convergence of pressure term $\Pi_{\delta}(\rho,n)$.
\begin{lem} \label{le4.1}
 Let  $ (\rho_\epsilon, n_\epsilon,\mathbf{u}_{\epsilon},\mathbf{H}_{\epsilon})$ be the sequence of solutions of the problem \eqref{app1}-\eqref{app42} constructed in Proposition \ref{prop2}, then there exists a constant $C$ independent of $\epsilon,$ such that		
\begin{equation}
\int_{0}^{T} \int_{\Omega}\left(n_\epsilon P \left(\rho_\epsilon,n_\epsilon\right)+\delta\left(n_\epsilon^{B+1}+n_\epsilon \rho_\epsilon^{B}\right)\right) \mathrm{d} x \leq C\left(T,\delta, \rho_{0}, n_{0},\mathbf{m}_{0},\mathbf{H}_{0}\right). \label{eest}
\end{equation}
\end{lem}
\noindent{\bf{Proof:}}
We can follow the same path as those in Section 3 of Feireisl et al. \cite{feireisl2001existence} to prove this lemma. Recall the Bogovskii operator
$$
\mathcal{B}=[\mathcal{B}_1,\mathcal{B}_2,\mathcal{B}_3]:\left\{f \in L^{p}(\Omega) \mid \int_{\Omega} f=0\right\} \mapsto\left[W_{0}^{1, p}(\Omega)\right]^{3}
$$
and its properties
$$\operatorname{div} \mathcal{B}(f)=f,~\|\mathcal{B}(f)\|_{W^{1, p}(\Omega)} \leq C\|f\|_{L^{p}(\Omega)},~\|\mathcal{B}(\operatorname{div} \mathbf{g})\|_{L^{q}(\Omega)} \leq C\|\mathbf{g}\|_{L^{q}(\Omega)}$$
for $1<p,q<\infty$. Let ${\bm \phi}$ be a test function on \eqref{app0}, where
${\bm \phi}=\left(\phi_1,\phi_2,\phi_3\right)$ and
$$\phi_i=\psi(t) \mathcal{B}_{i}\left[n_{\varepsilon}-m_{0}\right],~ \psi\in \mathcal{D}(0, T),~ 0 \leq \psi \leq 1, ~m_{0}=\frac{1}{|\Omega|} \int_\Omega n(t) \mathrm{d} x,~i=1,2,3.$$
 Noticing that
\begin{align}
	\int_{0}^{T} \int_{\Omega}& \Pi_{\delta}(\rho_{k}, n_{k}) \operatorname{div} {\bm\phi} d x d t=
\int_{0}^{T} \int_{\Omega} \psi(t) \operatorname{div} \mathcal{B}\left(n_{\epsilon}-m_{0}\right) \Pi_{\delta}(\rho,n) d x d t\nonumber\\
&=\int_{0}^{T} \int_{\Omega} \psi(t)\left( n_{\epsilon}\Pi_{\delta}(\rho_\epsilon, n_\epsilon) -m_{0}\Pi_{\delta}(\rho_\epsilon, n_\epsilon)\right) d x d t,\\
\int_0^T\int_{\Omega} \psi(t)&((\nabla \times \mathbf{H}_{\epsilon}) \times \mathbf{H}_{\epsilon}) \cdot \mathcal{B}\left[n_{\epsilon}-m_{0}\right] \mathrm{d} x \mathrm{d} t\nonumber\\
&=-\int_0^T\int_{\Omega} \psi(t)\left(\mathbf{H}_{\epsilon}^{\top} \nabla  \mathcal{B}\left[n_{\epsilon}-m_{0}\right] \mathbf{H}_{\epsilon}+\frac{1}{2} \nabla\left(|\mathbf{H}|_{\epsilon}^{2}\right) \cdot \mathcal{B}\left[n_{\epsilon}-M_{0}\right] \right) \mathrm{d} x \mathrm{d} t\nonumber\\
&=-\int_0^T\int_{\Omega} \psi(t)\mathbf{H}_{\epsilon}^{\top} \nabla  \mathcal{B}\left[n_{\epsilon}-m_{0}\right] \mathbf{H}_{\epsilon}\mathrm{d} x \mathrm{d} t+\int_0^T\int_{\Omega} \psi(t)\frac{1}{2}|\mathbf{H}|_{\epsilon}^{2} \cdot \left(n_{\epsilon}-m_{0} \right) \mathrm{d} x \mathrm{d} t,
\end{align}
and
\begin{align}
m_0\int_{0}^{T} \int_{\Omega} \psi(t)\Pi_{\delta}&(\rho_\epsilon, n_\epsilon)\mathrm{d} x \mathrm{d} t \leq C,\nonumber\\
\int_0^T\int_{\Omega} \psi(t)\mathbf{H}_{\epsilon}^{\top}& \nabla  \mathcal{B}\left[n_{\epsilon}-m_{0}\right] \mathbf{H}_{\epsilon}\mathrm{d} x \mathrm{d} t+\int_0^T\int_{\Omega} \psi(t)\frac{1}{2}|\mathbf{H}|_{\epsilon}^{2} \cdot n_{\epsilon}\mathrm{d} x \mathrm{d} t\nonumber\\
&\leq C \int_0^{T}\psi(t)\|\mathbf{H}_{\epsilon}\|_{L^4}^2\|n_{\epsilon}\|_{L^2} \mathrm{d} t\nonumber\\
&\leq C \|\nabla \mathbf{H}_{\epsilon}\|_{L^2L^2}+\|n_{\epsilon}\|_{L^\infty L^2}\nonumber\\
&\leq C\left(T,\delta, \rho_{0}, n_{0},\mathbf{m}_{0},\mathbf{H}_{0}\right),
\end{align}
we can obtain \eqref{eest} following the method in  \cite{feireisl2001existence}, as the rest is almost the same.
 \begin{flushright}
      $\square$
      \end{flushright}
With \eqref{appest2}, \eqref{appest3}, \eqref{eest} and letting  $\epsilon$ of the \eqref{app1}-\eqref{app41} pass to $0^+$, we know the limit  $(\rho,n,\mathbf{u},\mathbf{H})$ solves following system on $Q_T$ in distribution sense
\begin{align} 
	\partial_{t} \rho&+\operatorname{div}(\rho \mathbf{u}) = 0, \label{e1}\\
	\partial_{t} n&+\operatorname{div}(n \mathbf{u}) =0, \label{e2}\\
	\partial_{t}((\rho+n) \mathbf{u})+\operatorname{div}((\rho+n) \mathbf{u} \otimes \mathbf{u})&+\nabla \overline{\Pi}_{\delta}(\rho, n)=\mu \Delta \mathbf{u}\nonumber\\
	&+(\mu+\lambda) \nabla \operatorname{div} \mathbf{u}+(\nabla \times \mathbf{H}) \times \mathbf{H}, \label{e3}\\
	\mathbf{H}_{t}-\nabla \times(\mathbf{u} \times \mathbf{H})=-&\nabla \times(\nu \nabla \times \mathbf{H}),~~~
	\operatorname{div} \mathbf{H}=0 \label{e4},\\
	\left.(\rho, n,(\rho+n)\mathbf{u},\mathbf{H})\right|_{t=0}&=\left(\rho_{0, \delta}, n_{0, \delta}, \mathbf{m}_{0, \delta},\mathbf{H}_{0, \delta}\right) \text { on } \overline{\Omega}, \label{e5}\\
	\left(\frac{\partial \rho}{\partial\nu}\right.&, \left. \left.\frac{\partial n}{\partial\nu},\mathbf{u},\mathbf{H}\right)\right|_{\partial \Omega} =0. \label{e6}
\end{align}
The process of passing the limit about magnetic field is the same as the analysis in Faedo-Garlerkin approach, and the rest terms can be treated in the same  manner in Feireisl et al. \cite{feireisl2001existence}. Indeed, when $\epsilon \rightarrow 0^{+}$, one has
\begin{align*}
&\left(\rho_{\epsilon}, n_{\epsilon}\right) \rightarrow(\rho, n) \text { in } C\left(0, T ; L_{\text {weak}}^{B}(\Omega)\right) \text { and } n_{\epsilon} \rightarrow n \text{ weakly in } L^{B+1}\left(Q_{T}\right), \\
&\left(\epsilon \Delta \rho_{\epsilon}, \epsilon \Delta n_{\epsilon}\right) \rightarrow 0 \text {  in } L^{2}\left(0, T ; H^{-1}(\Omega)\right),~~
\mathbf{u}_{\epsilon} \rightarrow u \text { weakly in } L^{2}\left(0, T ; H_{0}^{1}(\Omega)\right),\\
&\left(\rho_{\epsilon}+n_{\epsilon}\right) \mathbf{u}_{\epsilon} \rightarrow(\rho+n) \mathbf{u} \text {  in } C\left([0, T] ; L_{\text {weak}}^{\frac{2 B}{B+1}}\right) \cap C\left([0, T] ; H^{-1}(\Omega)\right), \\
&\left(\rho_{\epsilon} \mathbf{u}_{\epsilon}, n_{\epsilon} \mathbf{u}_{\epsilon}\right) \rightarrow(\rho \mathbf{u}, n \mathbf{u}),~~\left(\rho_{\epsilon}+n_{\epsilon}\right) \mathbf{u}_{\epsilon} \otimes \mathbf{u}_{\epsilon} \rightarrow(\rho+n) \mathbf{u} \otimes \mathbf{u} \text {  in } \mathcal{D}^{\prime}\left(Q_{T}\right), \\
&\left(\nabla \times \mathbf{H}_{\epsilon}\right) \times \mathbf{H}_{\epsilon} \rightarrow(\nabla \times \mathbf{H}) \times \mathbf{H} ,~~\nabla \times\left(\mathbf{u}_{\epsilon} \times \mathbf{H}_{\epsilon}\right) \rightarrow \nabla \times(\mathbf{u} \times \mathbf{H}) \text { in } \mathcal{D}^{\prime}(\Omega \times(0, T)),\\
&\Pi_{\delta}(\rho_\epsilon, n_\epsilon)\rightarrow \overline{\Pi_{\delta}(\rho, n)  }\text { weakly in } L^{1}\left(Q_{T}\right),\\
&\epsilon \nabla u_{\epsilon} \cdot \nabla\left(\rho_{\epsilon}+n_{\epsilon}\right) \rightarrow 0 \text {  in } L^{1}\left(Q_{T}\right),
\end{align*}
where $\overline{\Pi_{\delta}(\rho, n)}$ denotes the weak limit of $\Pi_{\delta}(\rho_\epsilon, n_\epsilon)$ in at least $L^1(Q_T)$. We may also pass the limit in \eqref{appest2} and \eqref{appest3} to get
\begin{align}
	\frac{1}{c_{0}} n(x, t) &\leq \rho(x, t) \leq c_{0} n(x, t), \label{eest2}\\
	\int_{\Omega}\left(\rho^{2}+ n^{2}+(\rho+ n)|\mathbf{u}|^2 +|\mathbf{H}|^2\right.&\left.+2\mathcal{H}_{\delta}(\rho, n)\right)\mathrm{d}x+\nonumber\\
	\int_0^{T}\int_{\Omega}\left(\mu|\nabla \mathbf{u}|^{2}+(\lambda+\mu)\right.&\left.(\operatorname{div} \mathbf{u})^{2}+\nu|\nabla \times \mathbf{H}|^{2}\right) \mathrm{d} x\mathrm{d}t
	\nonumber
	\\& \leq E_\delta(0)+ C(\Sigma,\mu)\Sigma T \exp^{C(\Sigma,\mu)\Sigma  T}. \label{eest3}
\end{align}

We are in the position to show $ \overline{\Pi_{\delta}(\rho, n)}= \Pi_{\delta}(\rho, n)$.  However, $\Pi_{\delta}(\rho, n)$ involves two functions and one can't apply the theory of Feireisl. et al. \cite{feireisl2001existence} directly.  Novotn{\'y} et al. \cite{novotny2020weak} introduced a method that using \eqref{s} to transform $\rho$ into $n s$, $P(\rho_\epsilon,n_\epsilon)$ into $P(n_\epsilon s_\epsilon,n_\epsilon)$, and $\Pi_\delta(\rho_\epsilon,n_\epsilon)$ into $\Pi_\delta(n_\epsilon s_\epsilon,n_\epsilon)$. Thus, one can use Lemma \ref{le2} and follow the line of \cite{feireisl2001existence} to obtain the result, if we can prove following lemma, where $s$ is fixed function when concerning the process of  passing limit $\epsilon \rightarrow 0$.  The following results and proof are inspired by \cite{novotny2020weak}.  

\begin{lem} \label{le4.2}
	\begin{align}
	\lim _{\epsilon \rightarrow 0}\Pi_{\delta}(\rho_\epsilon , n_\epsilon )=\lim _{\epsilon \rightarrow 0}\Pi_{\delta}(n_\epsilon s, n_\epsilon ) \label{eap}
	\end{align}
holds on $Q_T$ in the weak sense.
	\end{lem}
\noindent{\bf{Proof:}}

The relationship $\rho_\epsilon=n_\epsilon s_\epsilon$  gives
\begin{align*}
	\Pi_{\delta}(\rho_\epsilon, n_\epsilon)=	\Pi_{\delta}(n_\epsilon s_\epsilon, n_\epsilon)=	\Pi_{\delta}(n_\epsilon s_\epsilon, n_\epsilon)-	\Pi_{\delta}(n_\epsilon s, n_\epsilon)+	\Pi_{\delta}(n_\epsilon s, n_\epsilon)
\end{align*}
and it follows mean value theorem and \eqref{Pd} that
\begin{align}
	&\lim _{\epsilon \rightarrow 0^+}  \int_{0}^{T} \int_{\Omega}\left|	\Pi_{\delta}(n_\epsilon s_\epsilon, n_\epsilon)-	\Pi_{\delta}(n_\epsilon s, n_\epsilon)\right| \mathrm{d} x \mathrm{d} t\nonumber\\
	\leq C(\delta)& \lim _{\epsilon \rightarrow 0^+}  \int_{0}^{T} \int_{\Omega}\left|	P(n_\epsilon s_\epsilon, n_\epsilon)-	P(n_\epsilon s, n_\epsilon)\right|+n_\epsilon^B\left|	s^B_\epsilon- s^B\right| \mathrm{d} x \mathrm{d} t\nonumber\\
	\leq C(\delta) \lim _{\epsilon \rightarrow 0^+} & \int_{0}^{T} \int_{\Omega}\left| \partial_{\rho}	P( \theta n_\epsilon s+(1-\theta)(n_\epsilon s_\epsilon-n_\epsilon s),n_\epsilon)\right|\left|s_\epsilon- s\right|n_\epsilon+n_\epsilon^B\left|	s^B_\epsilon- s^B\right| \mathrm{d} x \mathrm{d} t. \label{eapp1}
\end{align}
When $k_1$ large enough, one can obtain
$$(B-\frac{1}{k_1})\frac{k_1}{k_1-1}< B+1,$$
and
\begin{align}
	\lim _{\epsilon \rightarrow 0^+} & \int_{0}^{T} \int_{\Omega}n_\epsilon^B\left|	s^B_\epsilon- s^B\right| \mathrm{d} x \mathrm{d} t \nonumber\\
	& \leq C\lim _{\epsilon \rightarrow 0^+}   \left(	 \int_{0}^{T} \int_{\Omega}n_\epsilon^{(B-\frac{1}{k_1})\frac{k_1}{k_1-1}} \mathrm{d} x \mathrm{d} t \right)^{\frac{k_1-1}{k_1}}\left(	\int_{0}^{T} \int_{\Omega}n_\epsilon\left|	s^B_\epsilon- s^B\right|^{k_1} \mathrm{d} x \mathrm{d} t \right)^{\frac{1}{k_1}}\nonumber\\
	& \leq C \lim _{\epsilon \rightarrow 0^+}   \left(	\int_{0}^{T} \int_{\Omega}n_\epsilon\left|	s_\epsilon- s\right|^{k_1} \mathrm{d} x \mathrm{d} t \right)^{\frac{1}{k_1}}=0, \label{eap2}
\end{align}
where we used  \eqref{appt}, $s_\epsilon, s\leq c_0$ and
\begin{align*}
	\int_{0}^{T} \int_{\Omega}n_\epsilon\left|	s^B_\epsilon- s^B\right|^{k_1} \mathrm{d} x \mathrm{d} t &\leq B^{k_1}\int_{0}^{T} \int_{\Omega} n_\epsilon\max\{s_\epsilon,s\}^{B-1}\left|	s_\epsilon- s\right|^{k_1} \mathrm{d} x \mathrm{d} t \\
	&\leq C\int_{0}^{T} \int_{\Omega} n_\epsilon\left|	s_\epsilon- s\right|^{k_1} \mathrm{d} x \mathrm{d} t. 
\end{align*}
The \eqref{Pr} shows
\begin{align}
&	\left| \partial_{\rho}	P( \theta n_\epsilon s+(1-\theta)(n_\epsilon s_\epsilon-n_\epsilon s),n_\epsilon)\right|\left|s_\epsilon- s\right|n_\epsilon\nonumber\\
&	\leq C\left\{\left[\theta n_\epsilon s+(1-\theta)\left|n_\epsilon s_\epsilon-n_\epsilon s\right|\right]^{\gamma^+-1}+n_\epsilon^{\gamma^--\frac{\gamma^-}{\gamma^+}}\right\}\left|s_\epsilon- s\right|n_\epsilon\nonumber\\
&	\leq C\left\{n_\epsilon^{\gamma^+}+n_\epsilon^{\gamma^--\frac{\gamma^-}{\gamma^+}+1}\right\}\left|s_\epsilon- s\right| \label{eap3}
\end{align}
for $s_\epsilon, s\leq c_0$  and $\theta \in(0,1)$.
Set $B\geq \max\{\gamma^+,\gamma^--\frac{\gamma^-}{\gamma^+}+1\}$, deriving in the same way as \eqref{eap2}, we have
\begin{align}
\lim _{\epsilon \rightarrow 0^+}\int_{0}^{T} \int_{\Omega}\left|	P(n_\epsilon s_\epsilon, n_\epsilon)-	P(n_\epsilon s, n_\epsilon)\right| \mathrm{d} x \mathrm{d} t=0 ,
\end{align}
which combined \eqref{eapp1} and \eqref{eap2} gives \eqref{eap}.
 \begin{flushright}
	$\square$
\end{flushright} 

Therefore, we only need to prove 
\begin{prop} \label{prop4.2}
\begin{align}
	\overline{\overline{\Pi_{\delta}(ns, n)}}:=\lim _{\epsilon \rightarrow 0^+}\Pi_{\delta}(n_\epsilon s, n_\epsilon )=\Pi_{\delta}(n s, n) \label{eapp0}
\end{align}
a.e. on $Q_T$.
\end{prop}

To prove this proposition, we have to use the effective viscous flux $\Pi_{\delta}(n_\epsilon s, n_\epsilon )-(2 \mu+\lambda)  \operatorname{div} \mathbf{u}_\epsilon $ to get following identity.
\begin{prop} \label{prop4.3}
	We denote $\overline{\overline{\Pi_{\delta}(n s, n) n}}:=\lim\limits_{\epsilon \rightarrow 0^+}\Pi_{\delta}(n_\epsilon s, n_\epsilon ) n_\epsilon$, then the identity
\begin{equation}
\overline{\overline{\Pi_{\delta}(n s, n) n}}-(2 \mu+\lambda) \overline{n \operatorname{div} \mathbf{u}}=\overline{\overline{\Pi_{\delta}(ns , n)}} n-(2 \mu+\lambda) n \operatorname{div} \mathbf{u} \label{eapp}
\end{equation}
holds a.e. on $Q_T$.
\end{prop}

\noindent{\bf{Proof:}}

We can follow the path of Feireisl et al. \cite{feireisl2001existence} to prove this lemma. We can use test function
\begin{equation}
	{\bm\varphi}(t, x)=\psi(t)\left(\nabla \Delta^{-1}\left(n_{\epsilon} {\bm\phi(x)}\right)\right)(t, x), \quad \psi \in C_{c}^{1}(0, T), {\bm\phi} \in C_{c}^{1}(\Omega)
\end{equation}
and
\begin{equation}
	{\bm\varphi}(t, x)=\psi(t)\left(\nabla \Delta^{-1}\left(n {\bm\phi(x)}\right)\right)(t, x), \quad \psi \in C_{c}^{1}(0, T), {\bm\phi} \in C_{c}^{1}(\Omega)
\end{equation}
act on \eqref{app3}, \eqref{e3} respectively. Subtract the two identity and pass the $\epsilon$ to $0^+$. Noticing that $\rho=ns$, $\rho_\epsilon=n_\epsilon s_\epsilon$, we can  use analysis similar with \eqref{eapp1} to obtain
\begin{align*}
	\int_{0}^{T} \int_{\Omega}\psi(t){\bm\phi}(x) \overline{\Pi_{\delta}(\rho, n) n} d x d t=	\int_{0}^{T} \int_{\Omega}\psi(t){\bm\phi}(x)\overline{\overline{\Pi_{\delta}(n s, n) n}} d x d t, \\
		\int_{0}^{T} \int_{\Omega}\psi(t){\bm\phi}(x)\overline{\Pi_{\delta}(\rho, n) }n d x d t=	\int_{0}^{T} \int_{\Omega}\psi(t){\bm\phi}(x)\overline{\overline{\Pi_{\delta}(n s, n) }}n d x d t.
\end{align*}
For the rest terms, we can use  the classical Mikhlin multiplier theorem and the Div-Curl Lemma of compensated compactness to dispose. The terms including magnetic field can be treated similar with those in Lemma \ref{le4.1} and we will see  there is no essential difference with the analysis in \cite{feireisl2001existence}. For more details, we refer to \cite{feireisl2001existence} and Novotn{\'y} et al. \cite{novotny2020weak}. 
 \begin{flushright}
	$\square$
\end{flushright} 

\noindent{\bf{Proof of Proposition \ref{prop4.2}:}}

Following the line in \cite{feireisl2001existence}, we multiply \eqref{app2} with $b'_k(n_\epsilon)$, where
$$b_{k}(z)=\left\{\begin{array}{ll}
	z \log \left(z+\frac{1}{k}\right) & |z| \leq k \\
	(k+1) \log \left(k+1+\frac{1}{k}\right) & |z| \geq k+1
\end{array} \in C^{1}(\mathbb{R})\right.$$
integrate on $(0,\tau)\times\Omega$ and let $k\rightarrow \infty$, $\epsilon\rightarrow 0^+$, we have
\begin{equation}
	\forall \tau \in (0,T), \int_{\Omega} \overline{n\ln n}(\tau, \cdot) \mathrm{d} x-\int_{\Omega} n_{0} \ln n_{0} \mathrm{~d} x \leq \int_{0}^{\tau} \int_{\Omega} \overline{n \operatorname{div} \mathbf{u}} \mathrm{d} x \mathrm{~d} t. \label{eapp6}
\end{equation}
As $(n,\mathbf{u})$ is the renormalized solution of \eqref{e2}, we conclude from Lemma \ref{le2.1} that 
\begin{equation}
	\forall \tau \in (0,\tau), \int_{\Omega} n \ln n(\tau, \cdot) \mathrm{d} x-\int_{\Omega} n_{0} \ln n_{0} \mathrm{~d} x=\int_{0}^{\tau} \int_{\Omega} \overline{n \operatorname{div} \mathbf{u}} \mathrm{d} x \mathrm{~d} t
\end{equation}
which combined with \eqref{eapp6} and \eqref{eapp} gives
\begin{align}
	\int_{\Omega}(\overline{n \ln n}-n\ln n)(\tau, \cdot) \mathrm{d} x \leq \frac{1}{2 \mu+\lambda} \int_{0}^{\tau} \int_{\Omega}\left(\overline{\overline{\Pi_{\delta}(n s, n)}} n-\overline{\overline{\Pi_{\delta}(ns, n) n}}\right) \mathrm{d} x \mathrm{d} t.	\label{eapp2}
\end{align}

From \eqref{Pd} and \eqref{Prs}, we know that 
\begin{align*}
	\Pi_{\delta}(ns, n)=P(ns, n)+\delta\left(n^{B}(1+s^B)+\frac{1}{2} n^{B} (s^{B-2}+s^{2})\right),
\end{align*}
where $n^B$ is a convex function and $	P(ns, n)$ is non-decreasing by Lemma \ref{le2}. It follows that
\begin{align}
	\int_\Omega(\overline{n \ln n}-n \ln n)&(\tau, \cdot)\mathrm{d} x\leq\frac{1}{\lambda+2\mu}\int_0^T\int_\Omega \overline{\overline{P(ns ,n )}}n-\overline{\overline{P (ns, n)n}}\mathrm{d} x \mathrm{d} t\nonumber\\
	&+\frac{\delta}{\lambda+2\mu}\int_0^T\int_\Omega\left\{ (1+s^B)\left(\overline{n^{B}}n-\overline{n^{B}n}\right)\right.\nonumber\\
	&~~~~+\frac{1}{2}(s^{B-2}+s^{2} )\left.\left(\overline{n^{B} }n- \overline{n^{B} n}\right)\right\}\rho\mathrm{d} x \mathrm{d} t\leq 0, \label{eapp3}
\end{align}
 where we used Proposition \ref{prop2.2} and  
 \begin{align*}
 	\overline{n^{B} }n\leq\overline{n^{B} n},\text{ a.e. on }Q_T,~~\text{ for } n^B\text{ is convex}.
 \end{align*}
Since $n\ln n\leq \overline{n \ln n}$, \eqref{eapp3} shows
\begin{align*}
	\int_{\Omega}(\overline{n \ln n}-n \ln n)\mathrm{d} x =0.
\end{align*}
It allow us to have strong convergence of $n$. By \eqref{Pr} and \eqref{Pn}, $P(ns,n)$ and $\Pi_{\delta}(ns, n)$ are continuous function of $n$, thus  \eqref{eapp0} holds and
 there exists a global weak solution  $(\rho,n,\mathbf{u},\mathbf{H})$ to
\begin{align} 
	\partial_{t} \rho&+\operatorname{div}(\rho \mathbf{u}) = 0, \label{ad1}\\
	\partial_{t} n&+\operatorname{div}(n \mathbf{u}) =0, \label{ad2}\\
	\partial_{t}((\rho+n) \mathbf{u})+\operatorname{div}((\rho+n) \mathbf{u} \otimes \mathbf{u})&+\nabla \Pi_{\delta}(\rho, n)=\mu \Delta \mathbf{u}\nonumber\\
	&+(\mu+\lambda) \nabla \operatorname{div} \mathbf{u}+(\nabla \times \mathbf{H}) \times \mathbf{H}, \label{ad3}\\
	\mathbf{H}_{t}-\nabla \times(\mathbf{u} \times \mathbf{H})=-&\nabla \times(\nu \nabla \times \mathbf{H}),~~~
	\operatorname{div} \mathbf{H}=0 \label{ad4},\\
	\left.(\rho, n,(\rho+n)\mathbf{u},\mathbf{H})\right|_{t=0}&=\left(\rho_{0, \delta}, n_{0, \delta}, \mathbf{m}_{0, \delta},\mathbf{H}_{0, \delta}\right) \text { on } \overline{\Omega}, \label{ad5}\\
	\left(\frac{\partial \rho}{\partial\nu}\right.&, \left. \left.\frac{\partial n}{\partial\nu},\mathbf{u},\mathbf{H}\right)\right|_{\partial \Omega} =0. \label{ad6}
\end{align}
for any given $T>0, \delta>0$ and $B\geq A+2$.
Here  $(\rho,n,\mathbf{u},\mathbf{H})$  satisfy \eqref{ad3} in the sense of
\begin{align}
	\int_{\Omega}(\rho_{k}+&n_{k}) \mathbf{u}_{k}\cdot{\bm\phi} d x-\int_{\Omega} \mathbf{m}_{0, \delta}\cdot{\bm\phi}d x=	\int_{0}^{T} \int_{\Omega}\left((\rho_{k}+n_{k})(\mathbf{u}_{k} \otimes \mathbf{u}_{k}): \nabla{\bm\phi}+\Pi_{\delta}(\rho_{k}, n_{k}) \operatorname{div} {\bm\phi}\right) \mathrm{d} x \mathrm{d} t \nonumber\\
	&-\int_{0}^{T} \int_{\Omega}\left(\mu \nabla \mathbf{u}_{k}: \nabla{\bm\phi}+(\mu+\lambda) \operatorname{div} \mathbf{u}_{k} \operatorname{div} {\bm\phi}-(\nabla \times \mathbf{H}_{k}) \times \mathbf{H}_{k} )\cdot{\bm\phi}\right) \mathrm{d} x \mathrm{d} t. \label{app0'}
\end{align}
Take ${\bm\phi}=\mathbf{u}_{k}$ and follow the path of deriving \eqref{ener1}, we know the solutions satisfy
\begin{equation}
	E(t)+h_\delta(\rho,n)+\int_{0}^{t} \int_{\Omega}\left(\mu|D \mathbf{u}|^{2}+(\lambda+\mu)(\operatorname{divu})^{2}+\nu|\nabla \times \mathbf{H}|^{2}\right) \mathrm{d} x \mathrm{d} s \leq E(0)+h_\delta(\rho,n)|_{t=0}\label{appest5},
\end{equation}
where $E(t), E(0)$ and $h_\delta(\rho,n)$ are given in  \eqref{ener1},\eqref{ener0},\eqref{h} respectively.
 Moreover,  
\begin{align}
		\frac{1}{c_{0}} n(x, t) \leq \rho(x, t) \leq c_{0} n(x, t),\label{appest4}
\end{align}
follows \eqref{appest2}.

\section{The vanishing viscosity limit $\delta\rightarrow 0^{+}$} \label{s5.}
 ~~~~
This section shall recover the weak solution to \eqref{bi1}-\eqref{bi7} by passing the limit of $(\rho_\delta,n_\delta,\mathbf{u}_\delta,\mathbf{H}_\delta)$ as $\delta \rightarrow 0^+$. We can get following estimates of pressure term uniformly for $\delta$ following Feireisl et al. \cite{feireisl2001existence}.
\begin{lem} \label{le5.1}
	Let $(\rho_\delta,n_\delta,\mathbf{u}_\delta,\mathbf{H}_\delta)$ be the solution stated in Propostion \ref{prop4.2}, then 
\begin{align}
		\begin{aligned}
		\int_{0}^{T} \int_{\Omega}\left(
		n_\delta^{\gamma^-+\gamma^-_{{Bog}}}+\rho_{\delta}^{\gamma^+}n_\delta^{\gamma^-_{Bog}}
	+	\delta\left(n_\delta^{\gamma^-_{Bog}}\rho_{\delta}^{B}+ n_{\delta}^{B+\gamma^-_{Bog}}\right)\right) \mathrm{d} x \mathrm{~d} t \leq C 
	\end{aligned}\label{est2}
\end{align}
where 
 $ \gamma^-_{{Bog}}:= \min \{1,\frac{2}{3} \gamma^{-}-1, \frac{\gamma^-}{3}\}$, $ 1\leq\gamma^+,\frac 95 \leq\gamma^- $ and $C$ is a positive constant that is independent of $\delta$. 
\end{lem}

One can derive the result by using test function
$${\bm\phi}=\psi(t) \mathcal{B}\left(n_{\delta}^{\gamma^-_{{Bog}}}-\frac{1}{|\Omega|} \int_{\Omega} n_{\delta}^{\gamma^-_{{Bog}}}\right)$$ 
to get
	\begin{align}
		\int_{0}^{T} \int_{\Omega}\left(P(\rho_\delta,n_\delta)n_\delta^{\gamma^-_{Bog}}+\delta\left(n_\delta^{\gamma^-_{Bog}}\rho_{\delta}^{B}+ n_{\delta}^{B+\gamma^-_{Bog}}\right)\right) \mathrm{d} x \mathrm{~d} t \leq C
\end{align}
respectively, which combined with \eqref{P} and \eqref{appest4} gives \eqref{est2}. The process is similar with those in Lemma \ref{le4.1} and we refer the proof to Novotn{\'y} et al.  \cite{novotny2020weak},  Feireisl et al. \cite{feireisl2001existence} and Hu et al. \cite{hu2010global}.

With \eqref{appest4}, \eqref{appest5} and \eqref{est2}, we can let $\delta\rightarrow 0^{+}$. Following the analysis as in Section \ref{s4.} and Lemma \ref{le5.1}, we have
\begin{align*}
	&\left(\rho_{\delta}, n_{\delta}\right) \rightarrow(\rho, n) \text { in } C\left(0, T ; L_{\text {weak}}^{\max\{\gamma^+,\gamma^-\}}(\Omega)\right) \text { and } n_{\delta} \rightarrow n \text{ weakly in } L^{\gamma^-+\Theta}\left(Q_{T}\right), \\
	&\mathbf{u}_{\delta} \rightarrow u \text { weakly in } L^{2}\left(0, T ; H_{0}^{1}(\Omega)\right),~~~~P_{\delta}(\rho_\delta, n_\delta)\rightarrow \overline{P_{\delta}(\rho, n)  }\text { weakly in } L^{1}\left(Q_{T}\right),\\
	&\left(\rho_{\delta}+n_{\delta}\right) \mathbf{u}_{\delta} \rightarrow(\rho+n) \mathbf{u} \text {  in } C\left([0, T] ; L_{\text {weak}}^{\frac{2 \gamma^+}{\gamma^++1}}\right) \cap C\left([0, T] ; H^{-1}(\Omega)\right), \\
	&\left(\rho_{\delta} \mathbf{u}_{\delta}, n_{\delta} \mathbf{u}_{\delta}\right) \rightarrow(\rho \mathbf{u}, n \mathbf{u}),~~\left(\rho_{\delta}+n_{\delta}\right) \mathbf{u}_{\delta} \otimes \mathbf{u}_{\delta} \rightarrow(\rho+n) \mathbf{u} \otimes \mathbf{u} \text {  in } \mathcal{D}^{\prime}\left(Q_{T}\right), \\
	&\left(\nabla \times \mathbf{H}_{\delta}\right) \times \mathbf{H}_{\delta} \rightarrow(\nabla \times \mathbf{H}) \times \mathbf{H} ,~~\nabla \times\left(\mathbf{u}_{\delta} \times \mathbf{H}_{\delta}\right) \rightarrow \nabla \times(\mathbf{u} \times \mathbf{H}) \text { in } \mathcal{D}^{\prime}(\Omega \times(0, T)),
\end{align*}
and the limit $(\rho,n,\mathbf{u},\mathbf{H})$ solves
	\begin{align} 
	\partial_{t} \rho&+\operatorname{div}(\rho \mathbf{u}) = 0, \label{d1}\\
	\partial_{t} n&+\operatorname{div}(n \mathbf{u}) =0, \label{d2}\\
	\partial_{t}((\rho+n) \mathbf{u})+\operatorname{div}((\rho+n) \mathbf{u} \otimes \mathbf{u})&+\nabla \overline{P}(\rho, n)=\mu \Delta \mathbf{u}\nonumber\\
	&+(\mu+\lambda) \nabla \operatorname{div} \mathbf{u}+(\nabla \times \mathbf{H}) \times \mathbf{H}, \label{d3}\\
	\mathbf{H}_{t}-\nabla \times(\mathbf{u} \times \mathbf{H})=-&\nabla \times(\nu \nabla \times \mathbf{H}),~~~
	\operatorname{div} \mathbf{H}=0 \label{d4},\\
	\left.(\rho, n,(\rho+n)\mathbf{u},\mathbf{H})\right|_{t=0}&=\left(\rho_{0}, n_{0}, \mathbf{m}_{0},\mathbf{H}_{0}\right) \text { on } \overline{\Omega}, \label{d5}\\
	\left(\frac{\partial \rho}{\partial\nu}\right.&, \left. \left.\frac{\partial n}{\partial\nu},\mathbf{u},\mathbf{H}\right)\right|_{\partial \Omega} =0 \label{d6}
\end{align}
and satisfy 
	\begin{align}
	\frac{1}{c_{0}} n(x, t) &\leq \rho(x, t) \leq c_{0} n(x, t),\label{appest6}\\
		\int_{0}^{T} \int_{\Omega}\left(n^{\gamma^-+\gamma^-_{{Bog}}}+\right.&\left.\rho^{\gamma^+}n^{\gamma^-_{Bog}}\right) \mathrm{d} x \mathrm{~d} t \leq C\left(T,\rho_0,n_0,\mathbf{m}_0,\mathbf{H}_0\right), \label{appest7}\\
\int_{\Omega}\left(\frac{1}{2} (\rho+n) \mathbf{u}^{2}\right.&\left.+H_{P}(\rho, n)(t,\cdot)+\frac{1}{2}|\mathbf{H}|^{2}\right) \mathrm{d} x\nonumber\\
+\int_{0}^{t} \int_{\Omega}&\left(\mu|D \mathbf{u}|^{2}+(\lambda+\mu)(\operatorname{divu})^{2}+\nu|\nabla \times \mathbf{H}|^{2}\right) \mathrm{d} x \mathrm{d} s \leq E(0),  \label{appest8}
\end{align}
where $E(0)$ is given in \eqref{ener0}.

Next, we introduce following Lemma to prepare for the process of passing the limit.
\begin{lem}
 If $  \frac 95\leq \gamma^-,1\leq \gamma^+ $, then for  weak solution  $(\rho,n)$ satisfy \eqref{appest6}, \eqref{appest7}, \eqref{appest8}, we have 
	\begin{align}
		 \begin{aligned}\label{dest1}
	 &\int_{0}^{T} \int_{\Omega} n^{\gamma^--\frac{\gamma^-}{\gamma^+}+1} \mathrm{d} x \mathrm{~d} t  \\
	 &\leq 	C \left\{\int_{0}^{T}  \int_{\Omega}\left(n^{\gamma^-+\gamma^-_{{Bog}}-\Xi}+\rho^{\gamma^+}n^{\gamma^-_{Bog}-\Xi}\right) \mathrm{d} x \mathrm{~d} t \right\}^{\frac{1}{k_2}},  
	\end{aligned}
	\end{align}	
where 
$
	\Xi= \min\left\{\gamma^-_{Bog},\frac{\gamma^-}{\gamma^+},\frac{1}{10}\right\}>0
$ and $k_2>1$ is the H\"older relationship number that depends on $\Xi$.
\end{lem}
\noindent{\bf{Proof:}} The \eqref{appest6} shows for any constant $ \kappa>0 $, there is $ n^\kappa\leq C \rho^\kappa $.
Noticing that the $ C $ in \eqref{dest1} is depend on $ T $ and $ \Omega $ is bounded, then  by H\"older inequality and \eqref{appest6}, we know that the \eqref{dest1} can be reduced to prove
\begin{align} \label{dest1'}
\gamma^--\frac{\gamma^-}{\gamma^+}+1+\Xi\leq \Theta:=\max\{\gamma^++\gamma^-_{Bog},\gamma^-+\gamma^-_{Bog}\}
\end{align}
for $  \frac 95\leq \gamma^-,1\leq \gamma^+ $.

It is easy to verify that 
\begin{align*}
	\gamma^-_{Bog}=\left\{\begin{array}{l}
		\frac23	\gamma^--1 , \text { if } ~\gamma^-<3, \\
		1, \quad \quad~~~~ \text { if } ~ \gamma^-\geq3.
	\end{array}\right.
\end{align*}
If $ 3\leq\gamma^- $, then there is 
\begin{align*}
	\gamma^--\frac{\gamma^-}{\gamma^+}+1+\Xi\leq	\gamma^-+1 =\gamma^-+\gamma^-_{Bog}.
\end{align*}
If $ \gamma^-\geq \gamma^+ $, then we immediately have
\begin{align*}
	\gamma^--\frac{\gamma^-}{\gamma^+}+1+\Xi\leq \gamma^-+\Xi\leq \gamma^-+\gamma^-_{Bog}.
\end{align*}
If  $ \frac 95\leq \gamma^-<3$ and $  \gamma^-< \gamma^+  $, we only need to prove $ 	\gamma^--\frac{\gamma^-}{\gamma^+}+1+\frac{1}{10}\leq\gamma^++\gamma^-_{Bog} $ or
\begin{align} \label{dest10}
2+\frac{1}{10}\leq\gamma^+-\frac13	\gamma^-+\frac{\gamma^-}{\gamma^+}.
\end{align}
As $ \dfrac{2y}{3}+\dfrac{9}{5y}\geq 2\sqrt{\dfrac65}>2\dfrac{1}{10}$ for $\forall~ y>0 $, then one has
\begin{align*}
\gamma^+-\frac13	\gamma^-+\frac{\gamma^-}{\gamma^+}\geq\frac23	\gamma^++\frac{9}{5\gamma^+}>2\dfrac{1}{10},
\end{align*}
this is \eqref{dest10}. 

Combining all above, we have obtained \eqref{dest1'}.
 \begin{flushright}
	$\square$
\end{flushright}

\begin{lem} \label{le5.3}
	\begin{align}
		\lim _{\delta \rightarrow 0^+} P(\rho_\delta, n_\delta)=\lim _{\delta \rightarrow 0^+} 	P(n_\delta s, n_\delta ) \label{dap0}
\end{align}
holds in the weak sense.
\end{lem}

\noindent{\bf{Proof:}}

  We can arrive at the conclusion by using \eqref{appest6}, \eqref{appest7}, \eqref{dest1}, H\"older inequality, Proposition \ref{prop2.3} and following the path of Lemma \ref{le4.2}, if we notice the facts below.   The mean value theorem gives
\begin{align}
		\lim _{\delta \rightarrow 0^+}  \int_{0}^{T} \int_{\Omega}&\left|	P(\rho_\delta , n_\delta)- P(n_\delta s, n_\delta)\right| \mathrm{d} x \mathrm{d} t=
	\lim _{\delta \rightarrow 0^+}  \int_{0}^{T} \int_{\Omega}\left|	P(n_\delta s_\delta, n_\delta)- P(n_\delta s, n_\delta)\right| \mathrm{d} x \mathrm{d} t\nonumber\\
	\leq C & \lim _{\delta \rightarrow 0^+}  \int_{0}^{T} \int_{\Omega}\left| \partial_{\rho}	P( \theta  n_\delta s+(1-\theta)(n_\delta s_\delta-n_\delta s),n_\delta)\right|\left|s_\delta- s\right|n_\delta \mathrm{d} x \mathrm{d} t , \label{dap}
\end{align}
for some $ \theta \in (0,1)$. The \eqref{Pr} indicates
\begin{align}
	&	\left| \partial_{\rho}	P( \theta n_\delta s+(1-\theta)(n_\delta s_\delta-n_\delta s),n_\delta)\right|\left|s_\delta- s\right|n_\delta\nonumber\\
	&	\leq C\left\{\left[\theta n_\delta s+(1-\theta)\left|n_\delta s_\delta-n_\delta s\right|\right]^{\gamma^+-1}+n_\delta^{\gamma^--\frac{\gamma^-}{\gamma^+}}\right\}\left|s_\delta- s\right|n_\delta\nonumber\\
	&	\leq C\left\{n_\delta^{\gamma^+}+n_\delta^{\gamma^--\frac{\gamma^-}{\gamma^+}+1}\right\}\left|s_\delta- s\right| \label{dap1}
\end{align}
due to  $ \theta \in (0,1)$, $s_\delta, s\leq c_0$. The \eqref{dest1'} shows
\begin{align*}
	(\gamma^+-\frac{1}{k_3})\frac{k_3}{k_3-1}< \gamma^++\gamma^-_{Bog},\\
	(\gamma^--\frac{\gamma^-}{\gamma^+}+1-\frac{1}{k_4})\frac{k_4}{k_4-1}< \Theta,
\end{align*}
if $k_3, k_4$ large enough. The H\"older inequality and Proposition \ref{prop2.3} give 
\begin{align*}
		\lim _{\delta \rightarrow 0^+} & \int_{0}^{T} \int_{\Omega}n_\delta^{\gamma^+}\left|	s_\delta- s\right| \mathrm{d} x \mathrm{d} t \\
		& \leq C\lim _{\delta \rightarrow 0^+}   \left(	 \int_{0}^{T} \int_{\Omega}n_\delta^{(\gamma^+-\frac{1}{k_3})\frac{k_3}{k_3-1}} \mathrm{d} x \mathrm{d} t \right)^{\frac{k_3-1}{k_3}}\left(	\int_{0}^{T} \int_{\Omega}n_\delta\left|	s_\delta- s\right|^{k_3} \mathrm{d} x \mathrm{d} t \right)^{\frac{1}{k_3}}\\
		& \leq C \lim _{\delta \rightarrow 0^+}   \left(	\int_{0}^{T} \int_{\Omega}n_\delta\left|	s_\delta- s\right|^{k_3} \mathrm{d} x \mathrm{d} t \right)^{\frac{1}{k_3}}=0, 
\end{align*}
where 
\begin{align*}
	 \int_{0}^{T} \int_{\Omega}n_\delta^{(\gamma^+-\frac{1}{k_3})\frac{k_3}{k_3-1}} \mathrm{d} x \mathrm{d} t &\leq C \left(\int_{0}^{T} \int_{\Omega}n_\delta^{\gamma^++\gamma^-_{Bog}} \mathrm{d} x \mathrm{d} t\right)^{\frac{(\gamma^+-\frac{1}{k_3})\frac{k_3}{k_3-1}}{\gamma^++\gamma^-_{Bog}}}
	 \\
	 &\leq C \left(\int_{0}^{T} \int_{\Omega}\rho_\delta^{\gamma^+}n_\delta^{\gamma^-_{Bog}} \mathrm{d} x \mathrm{d} t\right)^{\frac{(\gamma^+-\frac{1}{k_3})\frac{k_3}{k_3-1}}{\gamma^++\gamma^-_{Bog}}}\leq C,
\end{align*}
due to \eqref{appest6}.
 \begin{flushright}
	$\square$
\end{flushright}  
We denote
\begin{align}
	\overline{\overline{P(ns, n)}}:=	\lim _{\delta \rightarrow 0^+}  P(n_\delta s, n_\delta )
\end{align}
and by Lemma \ref{le5.3}, we only need to prove
\begin{prop}
	\begin{align}
		\overline{\overline{P(n s, n)}}=P(n s, n)=P(\rho,n), \label{dapp0}
	\end{align} 
a.e. on $Q_T$.
\end{prop}

Recall the cut-off function family in \cite{feireisl2001existence}, 
$$
T_{k}(z)=k T\left(\frac{z}{k}\right), z \in \mathbb{R}, k=1,2, \cdots
$$
where $T \in C^{\infty}(\mathbb{R})$ satisfying
$$T(z)=\left\{\begin{array}{ll}
	z, & z \leq 1 \\
	2, & z \geq 3
\end{array} \in C^{\infty}(\mathbb{R}),\right. \text { concave, } z \in \mathbb{R}.$$
And, Lemma \ref{le2.1} suggest $(\rho_\delta, \mathbf{u}_\delta)$ is a renormalized solution of \eqref{ad1}. Thus, we have
\begin{equation}
	\partial_{t} T_{k}\left(n_{\delta}\right)+\operatorname{div}\left(T_{k}\left(n_{\delta}\right)\mathbf{u}_{\delta}\right)+\left[T_{k}^{\prime}\left(n_{\delta}\right) n_{\delta}-T_{k}\left(n_{\delta}\right)\right] \operatorname{div} \mathbf{u}_{\delta}=0 \quad \text { in } \mathscr{D}^{\prime}\left(Q_{T}\right)
\end{equation}
and let $\delta \rightarrow 0^+$, one yields
\begin{equation}
	\partial_{t} \overline{T_{k}(n)}+\operatorname{div}\left(\overline{T_{k}(n)} \mathbf{u}\right)+\overline{\left[T_{k}^{\prime}(n) n-T_{k}(n)\right] \operatorname{div} \mathbf{u}}=0 \text { in } \mathscr{D}^{\prime}\left(Q_{T}\right)
\end{equation}

And,  we will have the effective viscous flux identity as 
\begin{prop}
Identity
	\begin{equation}
		\overline{\overline{P(n s, n) T_{k}(n)}}-(2 \mu+\lambda) \overline{T_{k}(n) \operatorname{div} \mathbf{u}}=\overline{\overline{P(n s, n)}} ~\overline{T_{k}(n)}-(2 \mu+\lambda)  \overline{T_{k}(n)} \operatorname{div} \mathbf{u} \label{dapp}
	\end{equation}
	holds a.e. on $Q_T$.
\end{prop}
The proof is similar those in  Feireisl et al. \cite{feireisl2001existence} and the analysis in Proposition \ref{prop4.3}. See also \cite{novotny2020weak}.

As defined in \cite{feireisl2001existence}, we introduce 
$$
L_{k}(z)=\left\{\begin{array}{l}
	z \log z, \quad 0 \leq z \leq k \\
	z \log k+z \int_{k}^{z} \frac{T_{k}(s)}{s^{2}} d s, \quad z \geq k
\end{array}\right.
$$
satisfying
$$
L_{k}(z)=\beta_{k} z-2 k \text { for all } z \geq 3 k
$$
where
$$
\beta_{k}=\log k+\int_{k}^{3 k} \frac{T_{k}(s)}{s^{2}} d s+\frac{2}{3}.
$$
Let  $b_{k}(z):=L_{k}(z)-\beta_{k} z$ where $b_{k}^{\prime}(z)=0$ for all large $z,$ then
$$
b_{k}^{\prime}(z) z-b_{k}(z)=T_{k}(z).
$$
When $\gamma^-\geq \frac{9}{5}$, we have $\gamma^-+\gamma^-_{Bog}\geq2$, and by Lemma \ref{le2.1}, $(\rho_\delta,\mathbf{u}_\delta),(\rho,\mathbf{u})$ are renormalized solution of \eqref{ad1} and \eqref{d1} respectively. Thus,
\begin{align}
	\partial_{t} L_{k}\left(n_{\delta}\right)+\operatorname{div}\left(L_{k}\left(n_{\delta}\right) \mathbf{u}_{\delta}\right)+T_{k}\left(n_{\delta}\right) \operatorname{div} \mathbf{u}_{\delta}=0 \quad \text { in } \mathscr{D}^{\prime}\left((0, T) \times \Omega\right),\\
	\partial_{t} L_{k}(n)+\operatorname{div}\left(L_{k}(n) \mathbf{u}\right)+T_{k}(n) \operatorname{div} \mathbf{u}=0 \quad \text { in } \mathscr{D}^{\prime}\left((0, T) \times \Omega\right). 
\end{align}
Subtract these two equation, integrate on $(0,\tau)\times\Omega$. Then, letting $\delta \rightarrow 0$ and combining \eqref{dapp} and \eqref{Prs}, one has 
\begin{align}
	\begin{aligned} \label{dapp1}
		\int_{\Omega}\left(\overline{L_{k}(n)}-L_{k}(n)\right)&(\tau, \cdot) \mathrm{d} x=\int_{0}^{\tau} \int_{\Omega}\left(T_{k}(n)-\overline{T_{k}(n)}\right) \operatorname{div} \mathbf{u} \mathrm{d} x \mathrm{d} t \\
		&+\frac{1}{2 \mu+\lambda} \int_{0}^{\tau} \int_{\Omega}\left(\overline{\overline{P(ns, n)}}~ \overline{T_{k}(n)}-\overline{\overline{P(ns, n) T_{k}(n)}}\right) \mathrm{d} x \mathrm{d} t 
	\\
		&\quad \leq \int_{0}^{\tau} \int_{\Omega}\left(T_{k}(n)-\overline{T_{k}(n)}\right) \operatorname{div} \mathbf{u} \mathrm{d} x \mathrm{d} t
\end{aligned}
\end{align}
For the  term of the last line, we have
\begin{align} \label{dapp3}
	\left\|T_{k}\left(n_{\delta}\right)-T_{k}(n)\right\|_{L^{2}\left(Q_{T}\right)}\|\operatorname{div} \mathbf{u}\|_{L^{2}\left(Q_{T}\right)} \leq C \limsup _{\delta \rightarrow 0}\left\|T_{k}\left(n_{\delta}\right)-T_{k}(n)\right\|_{L^{1}\left(Q_{T}\right)}^{\frac{\gamma^--1}{2 \gamma^-}}=0.
\end{align}
Here, we used \eqref{appest7} and following lemma introduced in Novotn{\'y} et al.  \cite{novotny2020weak}. The proof is essentially given in \cite{novotny2020weak}, we write here for reader's convenience.
\begin{prop} \label{prop5}
 The sequence $n_{\delta}$ satisfies
\begin{align} \label{dapp00}
		\operatorname{osc}_{\gamma^-+1}\left[n_{\delta} \rightarrow n\right]\left(Q_{T}\right):=\sup _{k>1} \lim _{\delta \rightarrow 0} \sup \int_{Q_{T}}\left|T_{k}\left(n_{\delta}\right)-T_{k}n)\right|^{\gamma^-+1} \mathrm{~d} x \mathrm{~d} t<\infty.
\end{align}
\end{prop}

\noindent{\bf{Proof:}}

By \eqref{dapp} and \eqref{Prs}, we have	
\begin{align} \label{dapp01}
		\frac{\underline{q}_1}{2}\int_{0}^{T} \int_{\Omega}\left(\overline{n^{\gamma^-} T_{k}(n)}\right.&\left.-\overline{n^{\gamma^-}}~ \overline{T_{k}(n)}\right) \mathrm{d} x \mathrm{~d} t
		+\int_{0}^{T} \int_{\Omega}\left(\overline{\pi(n, s) T_{k}(n)}-\overline{\pi(n, s)}~ \overline{T_{k}(\rho)}\right) \mathrm{d} x \mathrm{~d} t\nonumber\\
	&	=
(2 \mu+\lambda)\lim _{\delta \rightarrow 0^+} \int_{0}^{T} \int_{\Omega}\left(T_{k}\left(\rho_{\delta}\right)-T_{k}(\rho)\right) \operatorname{div} \mathbf{u}_{\delta} \mathrm{d} x \mathrm{~d} t\\
&~~~+
		(2 \mu+\lambda) \lim _{\delta \rightarrow 0^+}\int_{0}^{T} \int_{\Omega}\left(T_{k}(\rho)-\overline{T_{k}(\rho)}\right) \operatorname{div} \mathbf{u}_{\delta} \mathrm{d} x \mathrm{~d} t\nonumber\\
		&=  \lim _{\delta \rightarrow 0^+} \left(I^1_\delta+I^2_\delta\right)\nonumber
\end{align}
Following the analysis in Feireisl et al. \cite{feireisl2001existence}, we have
\begin{align} \label{dapp02}
	\lim _{\delta \rightarrow 0^+} \left(I^1_\delta+I^2_\delta\right)&\leq 
	C \sup _{\delta>0}\left\|\operatorname{div} u_{\delta}\right\|_{L^{2}\left(Q_{T}\right) } \lim _{\delta \rightarrow 0}\left(\left\|T_{k}\left(n_{\delta}\right)-T_{k}(n)\right\|_{L^{2}\left(Q_{T}\right)}+\| T_{k}(n)-\overline{T_{k}(n)} \|_{L^{2}\left(Q_{T}\right)}\right) \nonumber\\
	&\leq 
	C \varlimsup_{\delta \rightarrow 0}\left\|T_{k}\left(n_{\delta}\right)-T_{k}(n)\right\|_{L^{2}\left(Q_{T}\right)}\nonumber\\
	&\leq 
	C \left[\operatorname{osc}_{\gamma^-+1}\left[n_{\delta} \rightarrow n\right]\left(Q_{T}\right)\right]^{\frac{1}{2\gamma^-}},
\end{align}
where we used Young inequality and Interpolation inequality in the last line.
Direct calculation shows
\begin{align} \label{dapp03}
		\int_{0}^{T} \int_{\Omega}&\left(\overline{n^{\gamma^-} T_{k}(n)}-\overline{n^{\gamma^-}}~ \overline{T_{k}(n)}\right) \mathrm{d} x \mathrm{~d} t\nonumber \\
		&=\limsup _{\delta \rightarrow 0} \int_{0}^{T} \int_{\Omega}\left(n_{\delta}^{\gamma^-}-n^{\gamma^-}\right)\left(T_{k}\left(n_{\delta}\right)-T_{k}(n)\right) \mathrm{d} x \mathrm{~d} t\nonumber \\
		&\quad+\int_{0}^{T} \int_{\Omega}\left(n^{\gamma^-}-\overline{n^{\gamma^-}}\right)\left(\overline{T_{k}(n)}-T_{k}(n)\right) \mathrm{d} x \mathrm{~d} t\nonumber \\
		&\geq \limsup _{\delta \rightarrow 0} \int_{0}^{T} \int_{\Omega}\left|T_{k}\left(n_{\delta}\right)-T_{k}(n)\right|^{\gamma^-+1} \mathrm{~d} x \mathrm{~d} t,
\end{align}
where we used convexity of $n \mapsto n^{\gamma^-}$ and concavity of $n \mapsto T_{k}(n)$ on $[0, \infty),$ and algebraic inequalities
$$
|a-b|^{\gamma} \leq\left|a^{\gamma}-b^{\gamma}\right| \text { and }|a-b| \geq\left|T_{k}(a)-T_{k}(b)\right|, \quad(a, b) \in[0, \infty)^{2}.
$$
With \eqref{dapp01}-\eqref{dapp03}, we have obtained \eqref{dapp00}.

\begin{flushright}
	$\square$
\end{flushright}  
As $\overline{L_{k}(n)} \rightarrow \overline{n \ln n}, L_{k}(n) \rightarrow n \ln n,$ in $C_{weak}\left([0, T] ; L^{r}(\Omega)\right)$ for any $1 \leq r<\gamma$, it follows \eqref{dapp1}-\eqref{dapp3} that
\begin{align}
	\int_{\Omega}(\overline{n \ln n}-n \ln n)(\tau, \cdot) \mathrm{d} x =0,\label{deapp6}
\end{align}
 which combined the  convexity of $n \ln n$  gives
\begin{align} \label{nd}
	n_{\delta} \rightarrow n \text { a.e. on } Q_{T}.
\end{align}
As $P(ns,n)$ is continuous function of $n$ due to \eqref{Pr} and \eqref{Pn}, thus  \eqref{dapp0} holds and the conclusion of Theorem \ref{thm1} has been obtained. Moreover, we can obtain $
\rho_{\delta} \rightarrow \rho \text { a.e. on } Q_{T}
$ by combining \eqref{appt} and \eqref{nd}.

\section{Proof of Theorem \ref{thm3}} \label{s6.}

~~~~~The method introduced in section \ref{s3.}-\ref{s5.} works for  $academic~magnetic~bi\text{-}fluid~ system$ \eqref{bi1}-\eqref{bi7} with pressure law \eqref{bi8} as well, and this  gives an approach that is different from those of Vasseur et al. \cite{vasseur2019global}. Moreover, the conclusion still holds when $\gamma\geq1,\alpha=\frac{9}{5}$, since we can control the amplitude of oscillations as Proposition \ref{prop5}. 

     Indeed, we can finish the Faedo-Galerkin approach by using the standard method in Feireisl et al. \cite{feireisl2001existence} and  Lemma \ref{le3.2}, as the explicit form of  pressure law  \eqref{bi8} won't produce essential difficulty  as in Section \ref{s3.}, see \cite{vasseur2019global}. For the process of passing to the limit of $\epsilon\rightarrow 0^{+}$ and $\delta\rightarrow 0^{+}$, we can use auxiliary function \eqref{s} 
     \begin{align*}
     		s:=\left\{\begin{array}{l}
     		\frac{\rho}{n} ~~\text {   if } n>0 \\
     		0~~~ \text {   if } n=0,
     	\end{array}\right.
     \end{align*}
 instead of 
 \begin{align*}
 	d=\rho+n, A=n/d, B=\rho/d
 \end{align*}
in \cite{vasseur2019global}. Let $P_\delta=P(\rho,n)+\delta(\rho+n)^\beta=\rho^\gamma+n^\alpha+\delta(\rho+n)^\beta$, then, 
\begin{align}
	P_\delta(\rho_\epsilon,n_\epsilon)=n_\epsilon^\gamma s_\epsilon^\gamma+n_\epsilon^\alpha+\delta(n_\epsilon s_\epsilon+n_\epsilon)^\beta,~~	P(\rho_\delta,n_\delta)=n_\delta^\gamma s_\delta^\gamma+n_\delta^\alpha \label{Ps}
\end{align}
can be treated similarly as Lemma \ref{le4.2} to arrive at 
	\begin{align}
			&\lim _{\epsilon \rightarrow 0^+}  P_\delta(n_\epsilon s_\epsilon, n_\epsilon)=	\lim _{\epsilon \rightarrow 0^+}  P_\delta(n_\epsilon s, n_\epsilon)\nonumber\\=\lim _{\epsilon \rightarrow 0^+} &(n_\epsilon^\gamma s^\gamma+n_\epsilon^\alpha)+\delta(n_\epsilon s+n_\epsilon)^\beta:=\overline{\overline{P_{\delta}(ns, n) }},\label{Ps1}\\
	\lim _{\delta \rightarrow 0^+}  P(n_\delta s_\delta, n_\delta )=&	\lim _{\delta \rightarrow 0^+}  P(n_\delta s, n_\delta )=\lim _{\delta \rightarrow 0^+} (n_\delta^\gamma s^\gamma+n_\delta^\alpha):=\overline{\overline{P(ns, n) }}, \label{Ps2}
\end{align} 
in the weak sense. Then, we can use the standard method in Feireisl et al. \cite{feireisl2001existence}  and part of the discussion of Section \ref{s4.}, \ref{s5.} to prove the strong convergence of $n$ and obtain the conclusion, if we have the proposition of controlling amplitude of oscillations. For instance, during the process of passing the limit $\epsilon\rightarrow 0^{+}$, we can prove
\begin{align*}
	\overline{\overline{P_{\delta}(ns, n) n}}-(2 \mu+\lambda)\overline{ n \operatorname{div} \mathbf{u}}:=\lim _{\epsilon \rightarrow 0^+} \left\{P(n_\epsilon s, n_\epsilon)n_\epsilon-(2 \mu+\lambda) n_\epsilon \operatorname{div} \mathbf{u}_\epsilon\right\}\\=\lim _{\epsilon \rightarrow 0^+} \left\{P(n_\epsilon s, n_\epsilon) -(2 \mu+\lambda)  \operatorname{div} \mathbf{u}_\epsilon \right\}\cdot\lim _{\epsilon \rightarrow 0^+} n_\epsilon:= \overline{\overline{P_{\delta}(ns, n) }}\cdot n-(2 \mu+\lambda)\overline{  \operatorname{div} \mathbf{u}}\cdot n
\end{align*}
	holds a.e. on $Q_T$, where $\overline{\overline{P_{\delta}(ns, n) n}}=	\lim\limits_{\epsilon \rightarrow 0^+} [n_\epsilon^\gamma s^\gamma+n_\epsilon^\alpha+\delta(n_\epsilon s+n_\epsilon)^\beta]n_\epsilon$. Then, noticing that $P_{\delta}(ns, n)$ is a non-deceasing function of $n$,  we can use similar discussion in Section \ref{s4.} to get
	\begin{align}
	0\leq\int_{\Omega}(\overline{n \ln n}-n\ln n)(\tau, \cdot) \mathrm{d} x \leq \frac{1}{2 \mu+\lambda} \int_{0}^{\tau} \int_{\Omega}\left(\overline{\overline{P_{\delta}(n s, n)}} n-\overline{\overline{P_{\delta}(ns, n) n}}\right) \mathrm{d} x \mathrm{d} t=0,
\end{align}
	which gives
	$$\overline{n \ln n}-n\ln n=0$$
	and strong convergence of $n$.
	
	Consequently, we only need to control amplitude of oscillations.
	\begin{prop} \label{prop6}
		The sequence $n_{\delta}$ satisfies
	\begin{align}
			\operatorname{osc}_{\alpha+1}\left[n_{\delta} \rightarrow n\right]\left(Q_{T}\right):=\sup _{k>1} \lim _{\delta \rightarrow 0} \sup \int_{Q_{T}}\left|T_{k}\left(n_{\delta}\right)-T_{k}n)\right|^{\alpha+1} \mathrm{~d} x \mathrm{~d} t<\infty. \label{o0}
	\end{align}

	\end{prop}
\noindent{\bf{Proof:}}

The Proposition \ref{prop2.2} shows
\begin{align} \label{o1}
	\begin{aligned}
	 &\int_{0}^{T} \int_{\Omega}\left(\overline{\overline{(ns)^\gamma T_k(n)}}-\overline{\overline{(ns)^\gamma}} ~\overline{T_k(n)}\right) \mathrm{d} x \mathrm{d} t\\
	 &= \int_{0}^{T} \int_{\Omega}\left(\lim _{\delta \rightarrow 0^+}n_\delta^\gamma s^\gamma T_k(n_\delta)-\lim _{\delta \rightarrow 0^+}n_\delta^\gamma s^\gamma\cdot\lim _{\delta \rightarrow 0^+}T_k(n_\delta)\right) \mathrm{d} x \mathrm{d} t \\
	 &\geq \int_{0}^{T} \int_{\Omega}\left(\lim _{\delta \rightarrow 0^+}n_\delta^\gamma T_k(n_\delta)-\lim _{\delta \rightarrow 0^+}n_\delta^\gamma\cdot\lim _{\delta \rightarrow 0^+}T_k(n_\delta)\right) s^\gamma \mathrm{d} x \mathrm{d} t\\
	 &=\int_{0}^{T} \int_{\Omega}\left(\overline{ T_k(n) n^\gamma }-  \overline{n^\gamma }~\overline{ T_k(n)} \right) s^\gamma \mathrm{d} x \mathrm{d} t\geq 0,
\end{aligned}
\end{align}
where $\overline{\overline{(ns)^\gamma}}=	\lim\limits_{\delta \rightarrow 0^+} n_\delta^\gamma s^\gamma$ and we used the fact that $z \mapsto T_{k}(z)$ and $z \mapsto z^{\gamma}$  are increasing functions in the last line of \eqref{o1}.
Define 
\begin{align*}
	\widetilde{H}_{\delta}:=P(n_\delta s, n_\delta)-(2 \mu+\lambda) \operatorname{div} \mathbf{u}_\delta,  \\
	\widetilde{H}=:\lim _{\delta \rightarrow 0^+} P(n_\delta s, n_\delta) -(2 \mu+\lambda)  \operatorname{div} \mathbf{u},
\end{align*}
where \begin{align}
	P(n_\delta s, n_\delta)=n_\delta^\gamma s^\gamma+n_\delta^\alpha.
\end{align}
 Then, we can prove
\begin{align}
\lim _{\delta \rightarrow 0^+} \int_{0}^{T} \psi \int_{\Omega} {\bm\phi} \widetilde{H}_{\delta}T_k(n_\delta)\mathrm{d} x \mathrm{d} t= \int_{0}^{T} \psi \int_{\Omega} {\bm\phi} \widetilde{H}\overline{T_k(n)}\mathrm{d} x \mathrm{d} t, \label{o2}
\end{align}
where $\overline{T_k(n)}=\lim\limits_{\delta \rightarrow 0^+}T_k(n_\delta)$, $\psi \in \mathcal{D}(0, T), {\bm\phi} \in D(\Omega)$.

The analysis above enable us to derive \eqref{o0} following the proof of Lemma 4.3 in Feireisl et al. \cite{feireisl2001existence} with slightly change. Indeed, as shown in Lemma 4.3 of \cite{feireisl2001existence}, one gets 
 $$\left|T_{k}(x)-T_{k}(y)\right|^{\alpha+1} \leq\left(x^{\alpha}-y^{\alpha}\right)\left(T_{k}(x)-T_{k}(y)\right), \quad \text{for}~\forall x, y \geq 0 .$$
 Then,
\begin{align}
	\begin{aligned}
	\lim_{\delta \rightarrow 0}& \|T_{k}\left(n_{\delta}\right)-T_{k}(n)\|_{L^{\alpha+1}\left(Q_{T}\right)}^{\alpha+1} \leq 
		\lim _{\delta \rightarrow 0} \int_{0}^{T} \psi \int_{\Omega} {\bm\phi}\left(n_{\delta}^{\alpha}-n^{\alpha}\right)\left(T_{k}\left(n_{\delta}\right)-T_{k}(n)\right) d x d t \\
		& \leq \lim_{\delta \rightarrow 0}\int_{0}^{T} \psi \int_{\Omega} {\bm\phi}\left(n_{\delta}^{\alpha}-n^{\alpha}\right)\left(T_{k}\left(n_{\delta}\right)-T_{k}(n)\right) d x d t \\
		&+\int_{0}^{T} \psi \int_{\Omega} {\bm\phi}\left(\overline{n^{\alpha}}-n^{\alpha}\right)\left(T_{k}(n)-\overline{T_{k}(n)}\right) d x d t\\
		&\leq \lim_{\delta \rightarrow 0}\int_{0}^{T} \psi \int_{\Omega} {\bm\phi}\left[n_{\delta}^{\alpha} T_{k}\left(n_{\delta}\right)-\overline{n^{\alpha}}~ \overline{T_{k}(n)}\right] d x d t\\
		&\leq\lim_{\delta \rightarrow 0} \int_{0}^{T} \psi \int_{\Omega} {\bm\phi}\left[\widetilde{H}_{\delta}+(2 \mu+\lambda)\operatorname{div} \mathbf{u}_{\delta}-n_\delta^\gamma s^\gamma\right] T_{k}\left(n_{\delta}\right)d x d t\\
		&\quad\quad-\int_{0}^{T} \psi \int_{\Omega} {\bm\phi}\left[\widetilde{H}+(2 \mu+\lambda) \operatorname{div} \mathbf{u} -\overline{\overline{(ns)^\gamma}}\right] \overline{T_{k}(n)} d x d t\\
		&\leq (2 \mu+\lambda)  \lim_{\delta \rightarrow 0}  \int_{0}^{T} \int_{\Omega}\left[\operatorname{div} \mathbf{u}_{\delta} T_{k}\left(n_{\delta}\right)-\operatorname{div} \mathbf{u} \overline{T_{k}(n)}\right] d x d t \\
		&\leq  C \lim _{\delta \rightarrow 0} \left| \int_{0}^{T} \int_{\Omega} \operatorname{div} \mathbf{u}_{\delta}\left[T_{k}\left(n_{\delta}\right)-T_{k}(n)+T_{k}(n)-\overline{T_{k}(n)}\right] d x d t\right| \\
		&\leq 
		 C \sup _{\delta>0}\left\|\operatorname{div} \mathbf{u}_{\delta}\right\|_{L^{2}\left(Q_{T}\right) } \lim _{\delta \rightarrow 0}\left(\left\|T_{k}\left(n_{\delta}\right)-T_{k}(n)\right\|_{L^{2}\left(Q_{T}\right)}+\| T_{k}(n)-\overline{T_{k}(n)} \|_{L^{2}\left(Q_{T}\right)}\right.  \\
			&\leq 
		 C \varlimsup_{\delta \rightarrow 0}\left\|T_{k}\left(n_{\delta}\right)-T_{k}(n)\right\|_{L^{2}\left(Q_{T}\right)}, \label{o3}
\end{aligned}
\end{align}
where we used $\overline{n^{\alpha}}\geq n^{\alpha}$, $T_{k}(n)\leq\overline{T_{k}(n)}$ in the third, fourth line and 
\eqref{o1}, \eqref{o2} in the fourth-seventh line. Finally, \eqref{o3} and Young inequality gives \eqref{o0}.

 \begin{flushright}
	$\square$
\end{flushright}  

\section*{Acknowledgement}
~~~ The authors would like to express their gratitude to Professor Milan Pokorn{\'y} for his kindness instructions and helpful discussions.

\begin{center}

\end{center}


\begin{thebibliography}{10}
	
	\bibitem{bresch2010global}
	D.~Bresch, B.~Desjardins, J.-M. Ghidaglia, and E.~Grenier.
	\newblock Global weak solutions to a generic two-fluid model.
	\newblock {\em Archive for Rational Mechanics and Analysis}, 196(2):599--629,
	2010.
	
	\bibitem{bresch2018multifluid}
	D.~Bresch, B.~Desjardins, J.-M. Ghidaglia, E.~Grenier, and M.~Hillairet.
	\newblock Multifluid models including compressible fluids.
	\newblock {\em Handbook of Mathematical Analysis in Mechanics of Viscous
		Fluids}, pages 2927--2978, 2018.
	
	\bibitem{bresch2012global}
	D.~Bresch, X.~Huang, and J.~Li.
	\newblock Global weak solutions to one-dimensional non-conservative viscous
	compressible two-phase system.
	\newblock {\em Communications in Mathematical Physics}, 309(3):737--755, 2012.
	
	\bibitem{bresch2018global}
	D.~Bresch and P.-E. Jabin.
	\newblock {Global existence of weak solutions for compressible Navier--Stokes
		equations: thermodynamically unstable pressure and anisotropic viscous stress
		tensor}.
	\newblock {\em Annals of Mathematics}, 188(2):577--684, 2018.
	
	\bibitem{bresch2019finite}
	D.~Bresch, P.~Mucha, and E.~Zatorska.
	\newblock {Finite-energy solutions for compressible two-fluid Stokes system}.
	\newblock {\em Archive for Rational Mechanics and Analysis}, 232(2):987--1029,
	2019.
	
	\bibitem{diperna1989ordinary}
	R.-J. DiPerna and P.-L. Lions.
	\newblock {Ordinary differential equations, transport theory and Sobolev
		spaces}.
	\newblock {\em Inventiones mathematicae}, 98(3):511--547, 1989.
	
	\bibitem{feireisl2004dynamics}
	E.~Feireisl.
	\newblock {\em Dynamics of viscous compressible fluids}, volume~26.
	\newblock Oxford University Press, 2004.
	
	\bibitem{feireisl2001existence}
	E.~Feireisl, A.~Novotn{\'y}, and H.~Petzeltov{\'a}.
	\newblock {On the existence of globally defined weak solutions to the
		Navier-Stokes equations}.
	\newblock {\em Journal of Mathematical Fluid Mechanics}, 3(4):358--392, 2001.
	
	\bibitem{hu2010global}
	X.~Hu and D.~Wang.
	\newblock {Global existence and large-time behavior of solutions to the
		three-dimensional equations of compressible Magnetohydrodynamic flows}.
	\newblock {\em Archive for Rational Mechanics and Analysis}, 197(1):203--238,
	2010.
	
	\bibitem{novotny2020weak}
	A.~Novotn{\'y} and M.~Pokorn{\'y}.
	\newblock Weak solutions for some compressible multicomponent fluid models.
	\newblock {\em Archive for Rational Mechanics and Analysis}, 235(1):355--403,
	2020.
	
	\bibitem{novotny2004introduction}
	A.~Novotn{\'y} and I.~Stra\v{s}kraba.
	\newblock {\em Introduction to the mathematical theory of compressible flow}.
	\newblock Number~27. Oxford University Press on Demand, 2004.
	
	\bibitem{vasseur2019global}
	A.~Vasseur, H.~Wen, and C.~Yu.
	\newblock Global weak solution to the viscous two-fluid model with finite
	energy.
	\newblock {\em Journal de Math{\'e}matiques Pures et Appliqu{\'e}es},
	125:247--282, 2019.
	
	\bibitem{wen2021global}
	H.~Wen.
	\newblock On global solutions to a viscous compressible two-fluid model with
	unconstrained transition to single-phase flow in three dimensions.
	\newblock {\em Calculus of Variations and Partial Differential Equations},
	60(4):1--38, 2021.
	
	\bibitem{huanyao2018review}
	H.~Wen, L.~Yao, and C.~Zhu.
	\newblock Review on mathematical analysis of some two-phase flow models.
	\newblock {\em Acta Mathematica Scientia}, 38(5):1617--1636, 2018.
	
	\bibitem{wen2018global}
	H.~Wen and L.~Zhu.
	\newblock Global well-posedness and decay estimates of strong solutions to a
	two-phase model with magnetic field.
	\newblock {\em Journal of Differential Equations}, 264(3):2377--2406, 2018.
	
\end{thebibliography}
\end{document}